\documentclass[letterpaper, 10pt,setspace]{elsarticle}

\usepackage[english]{babel}
\usepackage{amsmath, amssymb, amsthm, mathrsfs, bm, mathtools}
\usepackage{abstract, booktabs, caption, graphicx, psfrag, subfig, float, cases, cancel, chngpage, emptypage, color, enumerate}  
\usepackage{setspace}
\usepackage{fullpage}
\usepackage[figuresright]{rotating}
\usepackage{multirow}
\usepackage{algorithm}
\usepackage{algpseudocode}

\theoremstyle{plain}

\newcommand{\p}{\partial}

\graphicspath{images}

\begin{document}

\title{\bfseries Tweed and wireframe:
accelerated relaxation algorithms\\
for multigrid solution of elliptic PDEs\\
on stretched structured grids}

\author{Thomas Bewley\corref{cor1}}
\ead{bewley@eng.ucsd.edu}
\address{University of California San Diego, La Jolla, California, USA}

\author{Ali Mashayek\corref{cor1}}
\ead{mashayek@ic.ac.uk}
\address{Imperial College London, London, UK}

\author{Daniele Cavaglieri}
\ead{daniele.cavaglieri85@gmail.com}
\address{Tesla, Palo Alto, California, USA}

\author{Paolo Luchini}
\ead{luchini@unisa.it}
\address{University of Salerno, Fisciano, Italy}

\cortext[cor1]{Joint first and corresponding authors}

\maketitle

\begin{abstract}
\noindent Two new relaxation schemes are proposed for the smoothing step in the geometric multigrid solution of PDEs on 2D and 3D stretched structured grids.
The new schemes are characterized by efficient line relaxation on branched sets of lines of alternating colour,
where the lines are constructed to be everywhere orthogonal to the local direction of maximum grid clustering.
Tweed relaxation is best suited for grid clustering near the boundaries of the computational domain, whereas wireframe relaxation is best suited for grid clustering near the centre of the computational domain.
On strongly stretched grids of these types, multigrid leveraging these new smoothing schemes significantly outperforms multigrid based on other leading relaxation schemes, such as checkerboard and alternating-direction zebra relaxation, for the numerical solution of large linear systems arising from the discretization of elliptic PDEs.
\end{abstract}

\textbf{keywords:} Geometric multigrid; Gauss-Seidel; zebra relaxation; elliptic PDEs.

\section{Introduction}\label{sec:Intro}

Geometric multigrid methods are among the fastest techniques available for the numerical solution of the large linear systems arising from the high-resolution discretization of elliptic PDEs on structured grids~\cite{brandt2011,wesseling1995introduction,hackbusch2013multi,mccormick1987multigrid}, and can generally be implemented efficiently on parallel architectures~\cite{brandt1981,osti_897960}. Among applications of geometric multigrid methods are fluid dynamics and heat transfer~\cite{brandt2011,hortmann1990finite,wesseling2001geometric,esmaily2018scalable}, and heterogeneous and porous materials~\cite{boffy2012efficient,boffy2014multigrid,kanschat2017geometric}. 
In such applications, grid non-uniformity is at times inevitable, as a consequence of large separations of scales in the problem under consideration, arising from both multiscale physics, and heterogeneity in fluid properties, heat distribution, material properties, etc. A classic example is wall-bounded thermal convection, which requires a significant increase in resolution near the walls to resolve the turbulent thermal boundary layer\footnote{Multigrid methods have also been extended to handle certain non-elliptic PDEs (see~\cite{trottenberg2000} and the references therein).}.

The convergence rates of geometric multigrid methods depend largely on the interpolation and restriction
operators employed within. 
Appropriate relaxation schemes for the smoothing step are essential to accelerate the convergence of multigrid methods. While such methods have proven efficient on structured uniform grids, their efficiency often deteriorates dramatically once grid stretching becomes appreciable~\cite{naik1993improved,larsson2005conditional,zubair2007multigrid,bin2010geometric}. 
As reviewed in Esmaily $et$ $al.$ (2018;\cite{esmaily2018scalable}), various techniques have been proposed to increase the efficiency and stability of multigrid methods on nonuniform structured grids~\cite{chan2000robust}; such techniques include mapping the
non-uniform grid to a uniform coarse grid~\cite{esmaily2018scalable}, 
agglomeration of neighbouring elements~\cite{chan1998agglomeration,lallemand1992unstructured}, decomposition of stencil points~\cite{de1990matrix}, energy minimization~\cite{wan1999energy}, and incomplete Gauss elimination~\cite{reusken1996multigrid}.
This paper examines two new relaxation schemes that are particularly well suited for geometric multigrid methods on stretched grids.

The outline of this paper is as follows.
Section~\ref{sec:Theory} discusses the theoretical preliminaries required to describe our new methods.
Section~\ref{sec:Tweedwireframe} introduces the technical details regarding the implementation of tweed and wireframe relaxation for the solution of 2D and 3D elliptic PDEs.
Section~\ref{sec:Multigrid} shows how to implement tweed and wireframe relaxation in the smoothing step of a multigrid algorithm,
and determines convergence factors for different combinations of smoothing operators, restriction operators, and the number of smoothing steps applied at each step in the multigrid cycle.
Section~\ref{sec:Results} presents convergence results of multigrid with tweed and wireframe relaxation applied to the solution of the Poisson equation, which are shown in certain cases to compare quite favourably to the other available smoothing approaches.
Sections \ref{sec:3d} and \ref{sec:Complementarity} discuss the extension of this work to 3D.
Conclusions and future work are discussed in Section~\ref{sec:Conclusions}.

\section{Theoretical preliminaries}\label{sec:Theory}
A typical model problem used to evaluate multigrid performance is the solution of the elliptic 2D Poisson equation with Dirichlet boundary conditions, written here in its general (heterogenous, anisotropic) form:
\begin{subequations}\label{prob:CP}
\begin{alignat}{3}
	\frac{\p}{\p x_i} \Big( \sigma_{ij} \frac{\p u}{\p x_j} \Big) &= f(x,y) \quad &&\textrm{in}\ \Omega = [0,L_x] \times [0,L_y], \label{prob:CP;a} \\
	u &= g(x,y) &&\textrm{on}\  \p\Omega, \label{prob:CP;b}
\end{alignat}
\end{subequations}
where $x$ and $y$ are the spatial coordinates, $\Omega$ is the domain of interest, $\p\Omega$ is the boundary of $\Omega$, and the matrix with components $\sigma_{ij}(x,y)$ is symmetric positive definite.
In isotropic media, $\sigma_{ij}(x,y)=c(x,y)\,\delta_{ij}$.
In homogeneous media, $\sigma_{ij}(x,y)$ is constant in $x$ and $y$.
A system of this form may be isotropic, homogeneous, both, or neither.

To demonstrate our new method, we will consider the discretization of \eqref{prob:CP} in the homogeneous isotropic case, with $\sigma_{ij}(x,y)=\delta_{ij}$, on a stretched rectilinear grid in which all cells are rectangles (in 2D) or rectangular cuboids (in 3D); the extension of this method to anisotropic, inhomogeneous systems, curvilinear grids, other elliptic PDEs, and other boundary conditions follows using standard methods.
Discretization of this problem on an $(n_x+1) \times (n_y+1)$ stretched grid using a second-order central finite difference method, with $x_i$ and $y_j$ denoting the grid coordinates in the $x$ and $y$ directions, respectively, and $u_{i,j}$ denoting the discretized value of $u$ at the $\{i,j\}$ gridpoint, leads to a five-point discretization of the Laplacian:
\begin{subequations}\label{prob:DP}
\begin{align}
	& W_i\,u_{i-1,j} + E_i\,u_{i+1,j} + S_j\,u_{j-1,j} + N_j\,u_{i,j+1} + C_{i,j}\,u_{i,j} = f_{i,j}, \quad  i = 2,\dots,n_x,\ \ j=2,\dots,n_y, \\
	& u_{1,*},\  u_{n_x+1,*},\  u_{*,1},\  u_{*,n_y+1}\  \textrm{specified},
\end{align}
\end{subequations}
where the {W}est, {E}ast, {S}outh, and {N}orth coefficients are
\begin{equation*}
	W_i = \frac{1}{\Delta x_i \Delta x_{i-1/2}},       \quad
	E_i = \frac{1}{\Delta x_i \Delta x_{i+1/2}}, \quad
	S_j = \frac{1}{\Delta y_j \Delta y_{j-1/2}}, \quad
	N_j = \frac{1}{\Delta y_j \Delta y_{j+1/2}},
\end{equation*}
and where the {C}enter coefficient is $C_{i,\,j} = - (W_i + E_i + S_j + N_j)$, with
\begin{align}
    \Delta x_{i-1/2} &= x_i - x_{i-1}, \quad \Delta x_{i}=(\Delta x_{i+1/2}+\Delta x_{i-1/2})/2= (x_{i+1} - x_{i-1})/2,\\
    \Delta y_{j-1/2} &= y_j - y_{j-1}, \quad \Delta y_{j}=(\Delta y_{j+1/2}+\Delta y_{j-1/2})/2= (y_{j+1} - y_{j-1})/2.
\end{align}
It is thus seen that, even if the PDE is homogeneous and isotropic, grid stretching causes the discretized Poisson equation \eqref{prob:DP} to be inhomogeneous (with coefficients varying as a function of position)
and anisotropic (with coefficients varying as a function of direction).  

\begin{figure} 
	\centering
	\includegraphics[height=0.65\columnwidth]{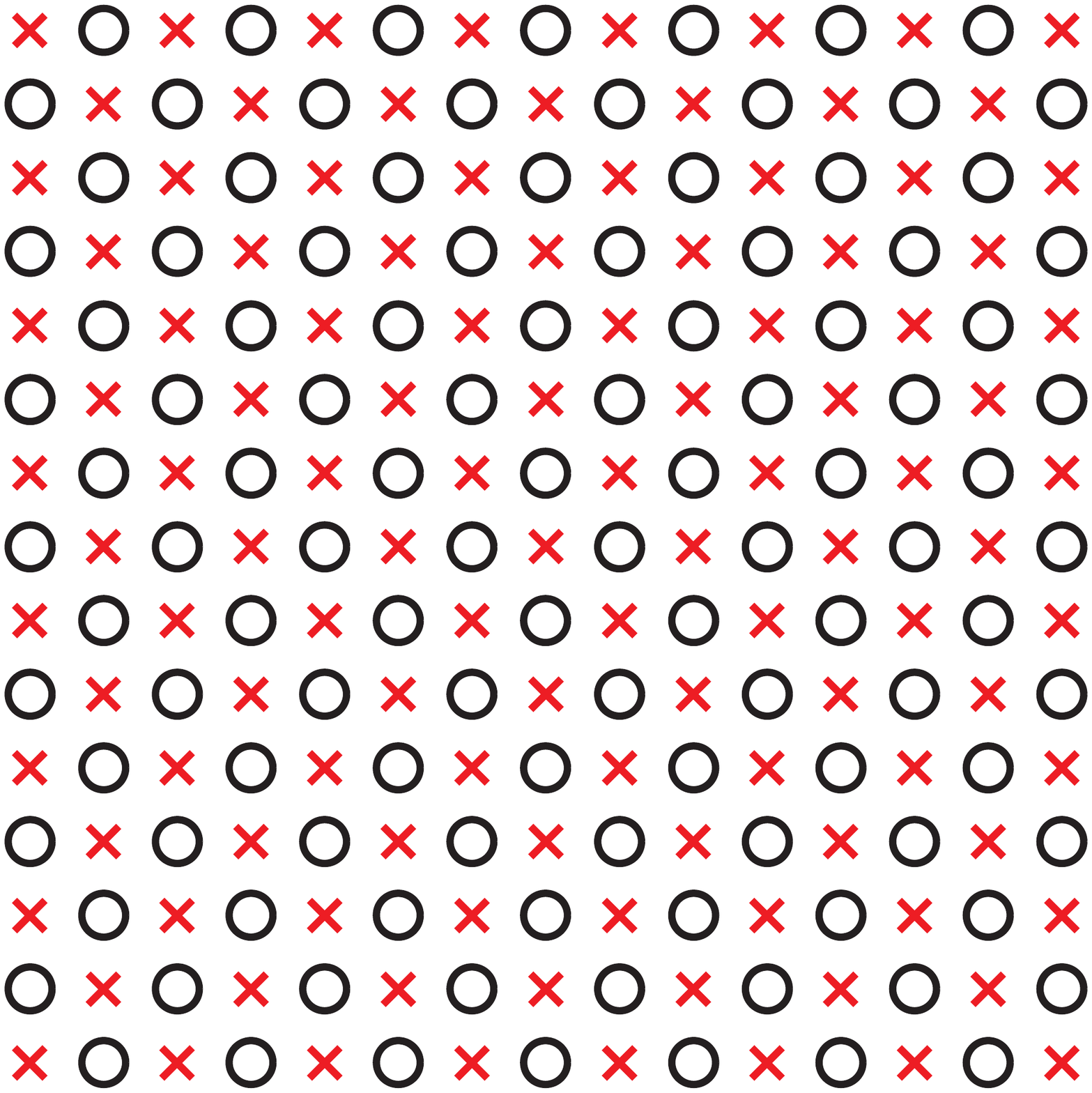}
	\caption{Gridpoint arrangement for checkerboard smoothing.}
	\label{Fig1}
\end{figure}

In the case of an unstretched grid (with $\Delta x$ and $\Delta y$ constant in both $x$ and $y$), we have
$W_i = E_i = {1}/{(\Delta x)^2}$ and $S_j = N_j = {1}/{(\Delta y)^2}$;
that is, the discretization of the Laplacian becomes homogeneous.  Further, if $\Delta x= \Delta y$, the discretization of the Laplacian also becomes isotropic, with $W_i = E_i = S_j = N_j$.
As is well known (see \cite{trottenberg2000}), for the case with no grid stretching,
standard smoothing approaches such as {\it checkerboard} (also called red-black point Gauss-Seidel) relaxation (see Figure~\ref{Fig1}) performs exceptionally well when applied within the multigrid framework
(see, e.g.,~\cite{luchini1994}).
As shown in \S \ref{sec:Results}, however, the performance of checkerboard relaxation starts to reduce significantly when grid stretching is introduced.

\begin{figure}  
	\centering
	\subfloat[\emph{Hyperbolic tangent stretching~\eqref{eq:tanh}, $c = 1.5$}]{\psfrag{x}[t][b]{\small$\tilde{x}$}\psfrag{y}[b][c]{\small$\tilde{y}$}\includegraphics[height=0.35\columnwidth]{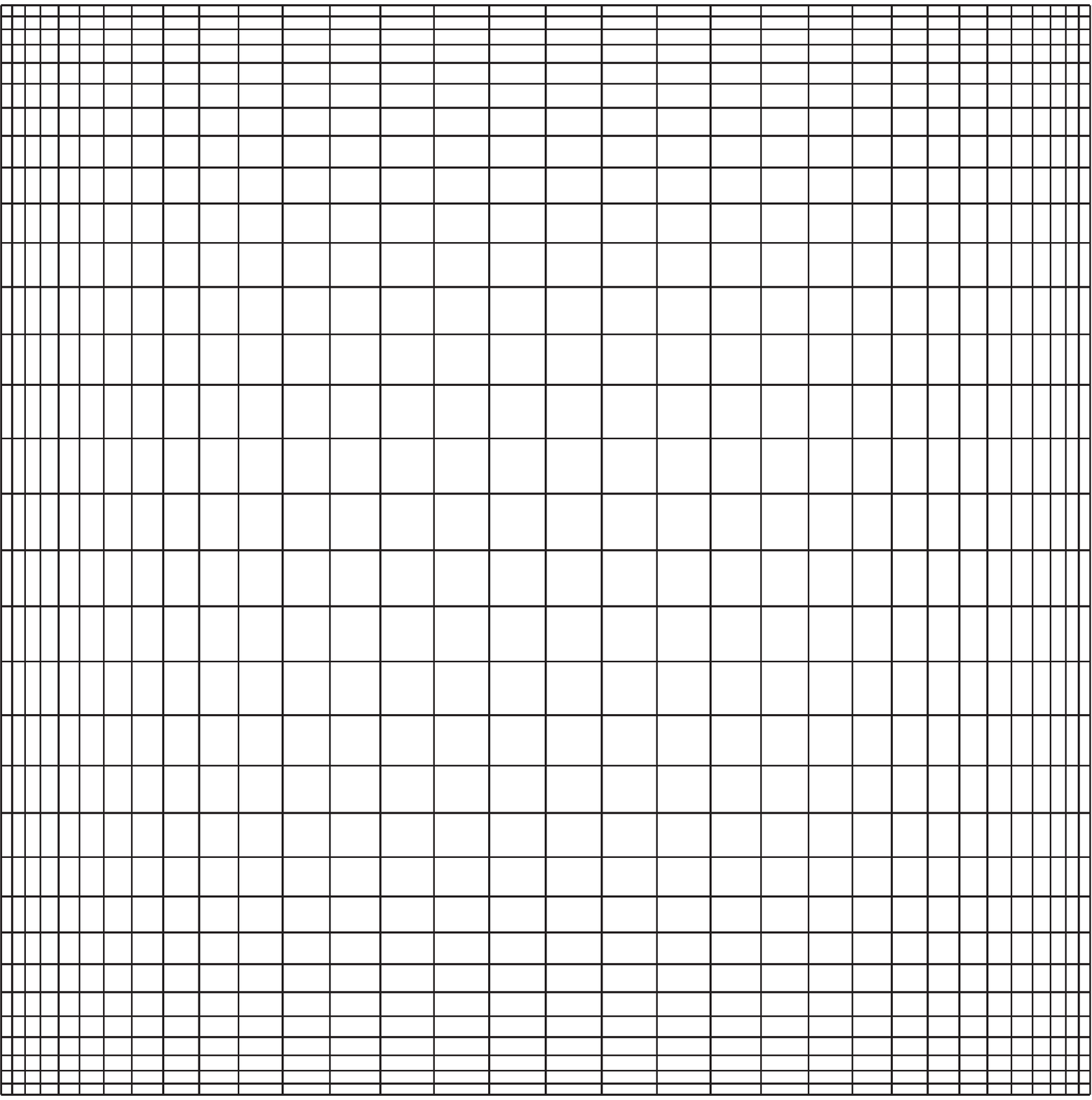}\label{fig:Str_grid_1p5}} \qquad\quad
	\subfloat[\emph{Hyperbolic tangent stretching~\eqref{eq:tanh}, $c = 3.0$}]{\psfrag{x}[t][b]{\small$\tilde{x}$}\psfrag{y}[b][c]{\small$\tilde{y}$}\includegraphics[height=0.35\columnwidth]{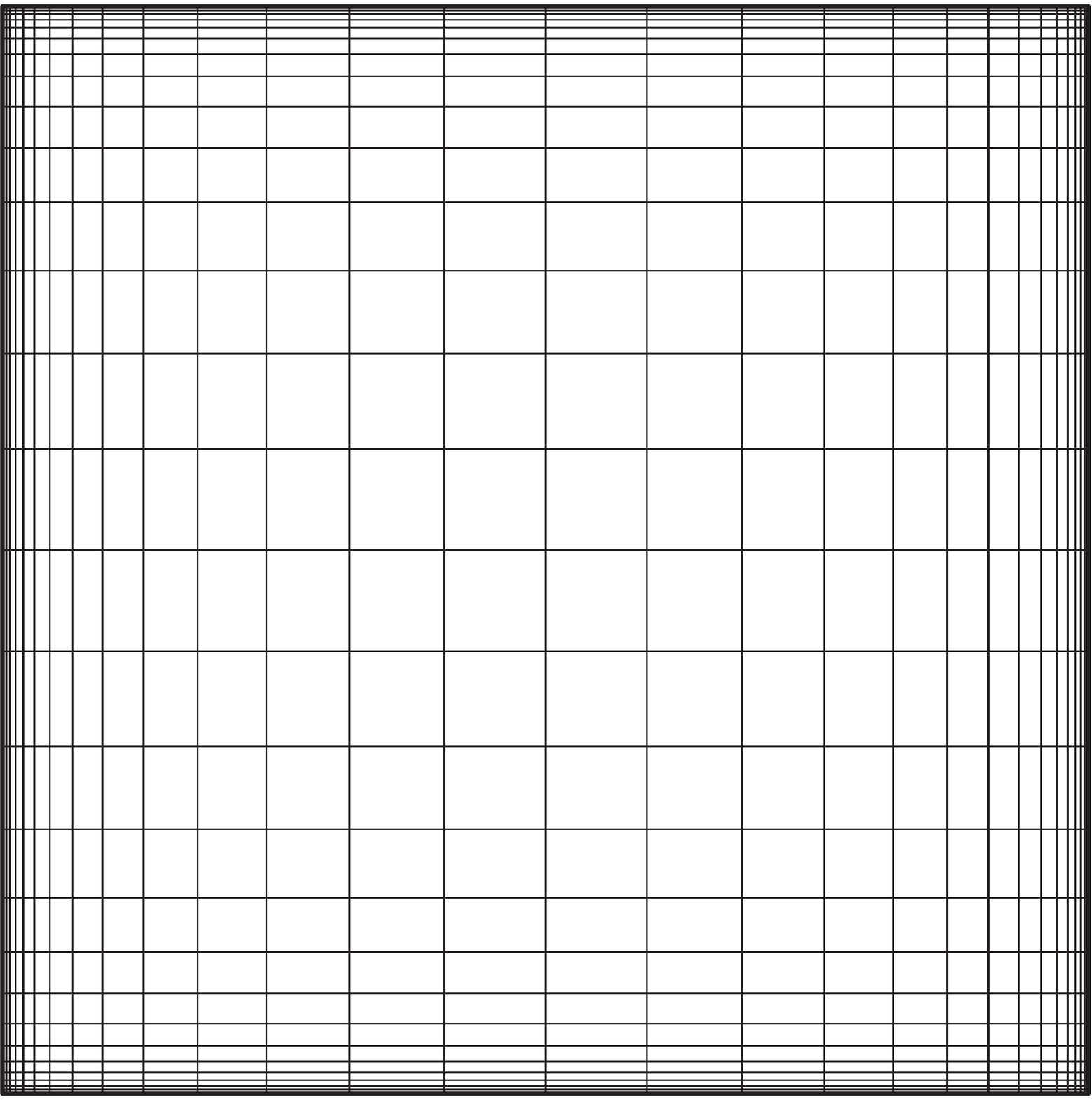}\label{fig:Str_grid_3p0}} \qquad\quad
	\subfloat[\emph{Hyperbolic tangent stretching~\eqref{eq:tanh_inv}, $c = 1.5$}]{\psfrag{x}[t][b]{\small$\tilde{x}$}\psfrag{y}[b][c]{\small$\tilde{y}$}\includegraphics[height=0.35\columnwidth]{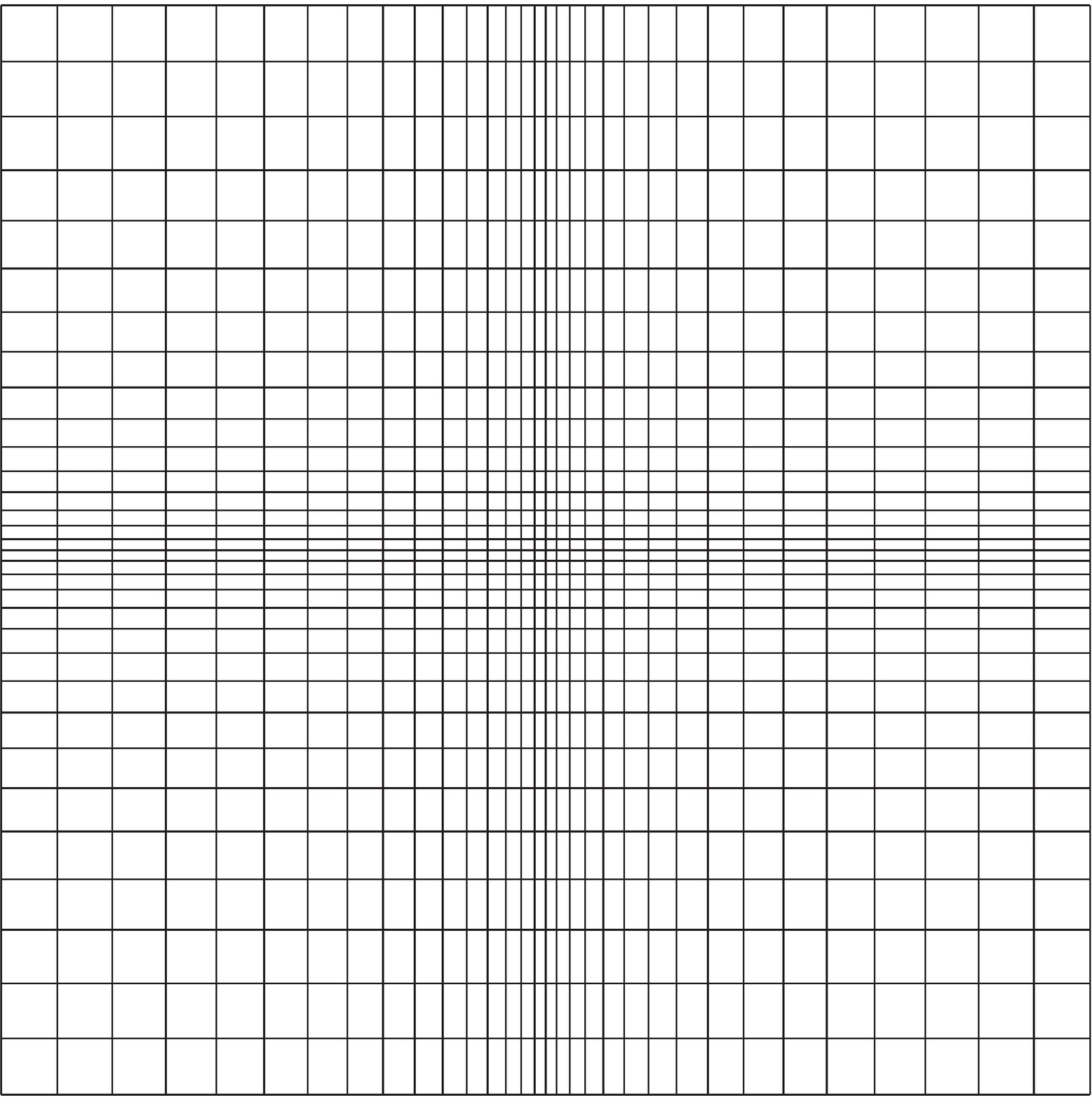}\label{fig:InvStr_grid_1p5}}
	\caption{2D stretched grids with clustering near the walls (a) and (b), and clustering near the centre (c).}
	\label{Fig2}
\end{figure}

When grid stretching is performed, gridpoints are clustered (that is, denser) in some directions more than others in certain regions of the computational domain, and thus two-delta waves in the error of the solution (that is, gridpoint-to-gridpoint oscillations, which checkerboard smoothing proves most effective at eliminating) are at mismatched spatial wavenumbers in different spatial directions in these regions.  When this happens, common wisdom~\cite{trottenberg2000} is that {\it zebra relaxation} (that is, line relaxation on lines of alternating ``colour'') should be applied in the direction orthogonal to the local direction of densest grid clustering. For example, in a square computational domain with grid clustering near the upper and lower boundaries, zebra relaxation along lines which are orthogonal to these two boundaries proves to be quite effective, whereas zebra relaxation along lines in the opposite direction proves to be much less effective.
In a 2D or 3D computational domain with grid clustering in {\it multiple} directions (see, e.g., Figure \ref{Fig2}), zebra relaxation in one direction alone is ineffective.
In such case, {\it alternating-direction zebra relaxation}, in which zebra relaxation is performed in each coordinate direction in succession (see Figure~\ref{Fig3}), is the method most commonly used,
and with it the rapid convergence of the multigrid approach may be recovered.
Note, however, that with the alternating-direction zebra relaxation approach, in any given region, the lines upon which relaxation is performed are orthogonal to the local direction of densest grid clustering
during only half of the sweeps in the 2D case, and during only a third of the sweeps in the 3D case.

Thus, rather than requiring two sweeps of zebra relaxation in 2D, or three sweeps of zebra relaxation in 3D, simply to get the relaxation lines used to be locally oriented to the local direction of densest grid clustering
during a fraction of the sweeps, we instead suggest a more effective motif for the line smoother.
The present work represents, we believe, the first attempt at developing a line smoother for multigrid relaxation which relaxes efficiently along {\it branched} lines,
that are locally orthogonal to the local directions of densest grid clustering
{\it everywhere} in the computational domain, even when grid clustering is applied in multiple directions, as indicated in Figure~\ref{Fig2} for the 2D case.

Towards this end, two different smoothers have been developed, which we denote \emph{tweed} and \emph{wireframe} relaxation.
Tweed relaxation efficiently addresses the problem of near-wall grid clustering (see Figure \ref{Fig2}a-b), whereas wireframe relaxation addresses grid clustering near the centre (see Figure \ref{Fig2}c).
The key idea behind the tweed motif (see Figure \ref{Fig4}) is to perform relaxation in blocks of connecting lines arranged in such a way as to make such lines {\it everywhere perpendicular to the closest domain boundary} if the domain is rectangular.
Red-black alternation, applied in a zebra-like fashion, makes the relaxation of each block independent from the other blocks of the same colour.
The key idea behind the wireframe motif (see Figure \ref{Fig5}), in contrast, is to perform relaxation in blocks of connecting lines arranged in such a way as to make such lines {\it everywhere perpendicular to the closest coordinate plane through the centre of the domain}
(and which are, thus, everywhere {\it parallel} to the closest domain boundary if the domain is rectangular).
Again, red-black alternation, applied in a zebra-like fashion, makes the relaxation of each block independent of the other blocks of the same colour.
Note that the tweed and wireframe motifs extend naturally to 3D, as discussed in section~\ref{sec:3d}.

\section{The new motifs: tweed and wireframe relaxation}\label{sec:Tweedwireframe}
\begin{figure}  
	\centering
	\subfloat[\emph{$x$-direction smoothing}]{\includegraphics[height=0.35\columnwidth]{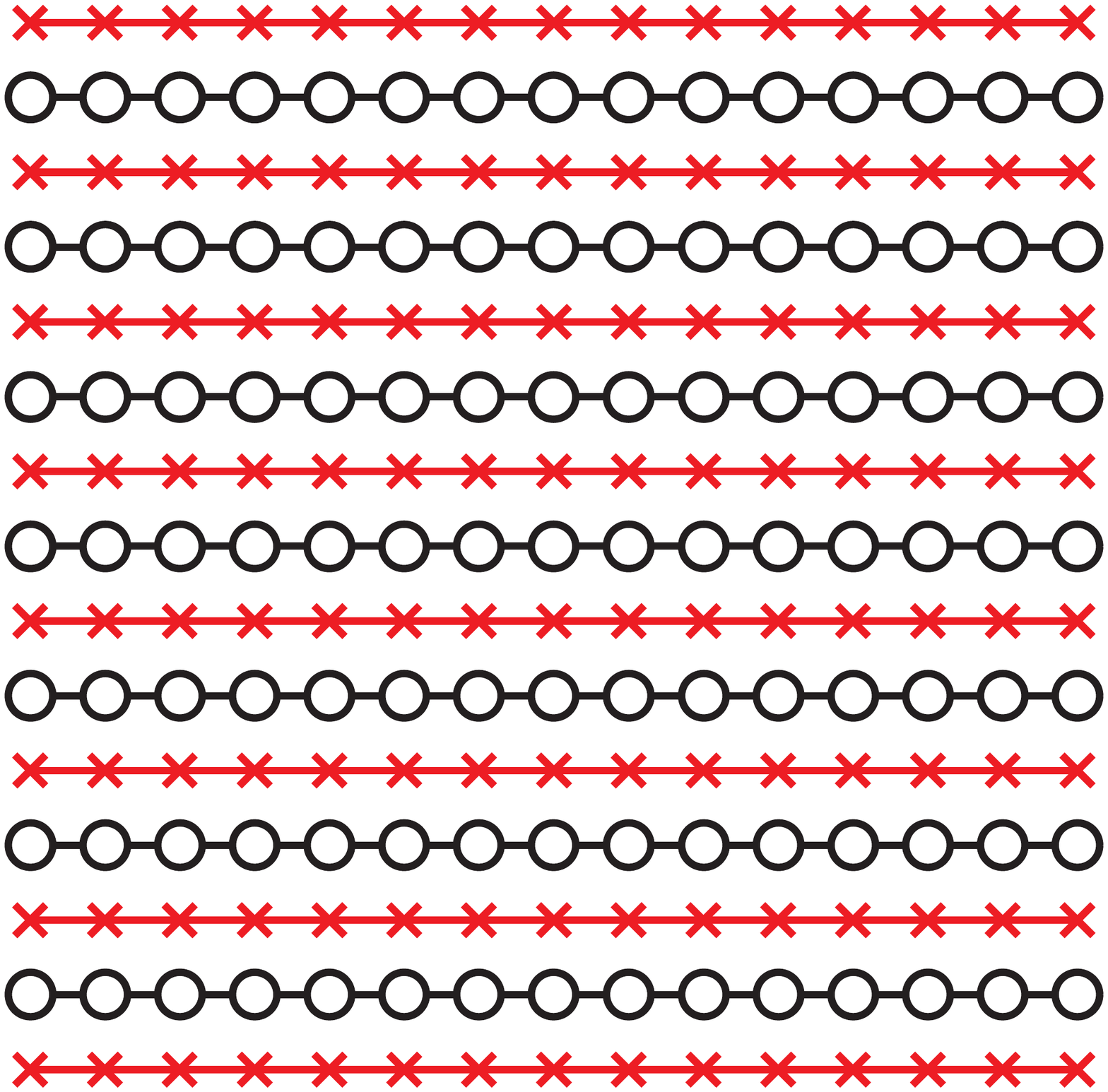}} \quad \subfloat[][\emph{$y$-direction smoothing}]{\includegraphics[height=0.35\columnwidth, angle=90, origin=c]{Figs/Fig3/zebrax}}
	\caption{Relaxation motifs for 2D alternating-direction zebra smoothing.}\label{Fig3}
	\subfloat[\emph{Square grid ($n_x = n_y$)}]{\includegraphics[height=0.35\columnwidth]{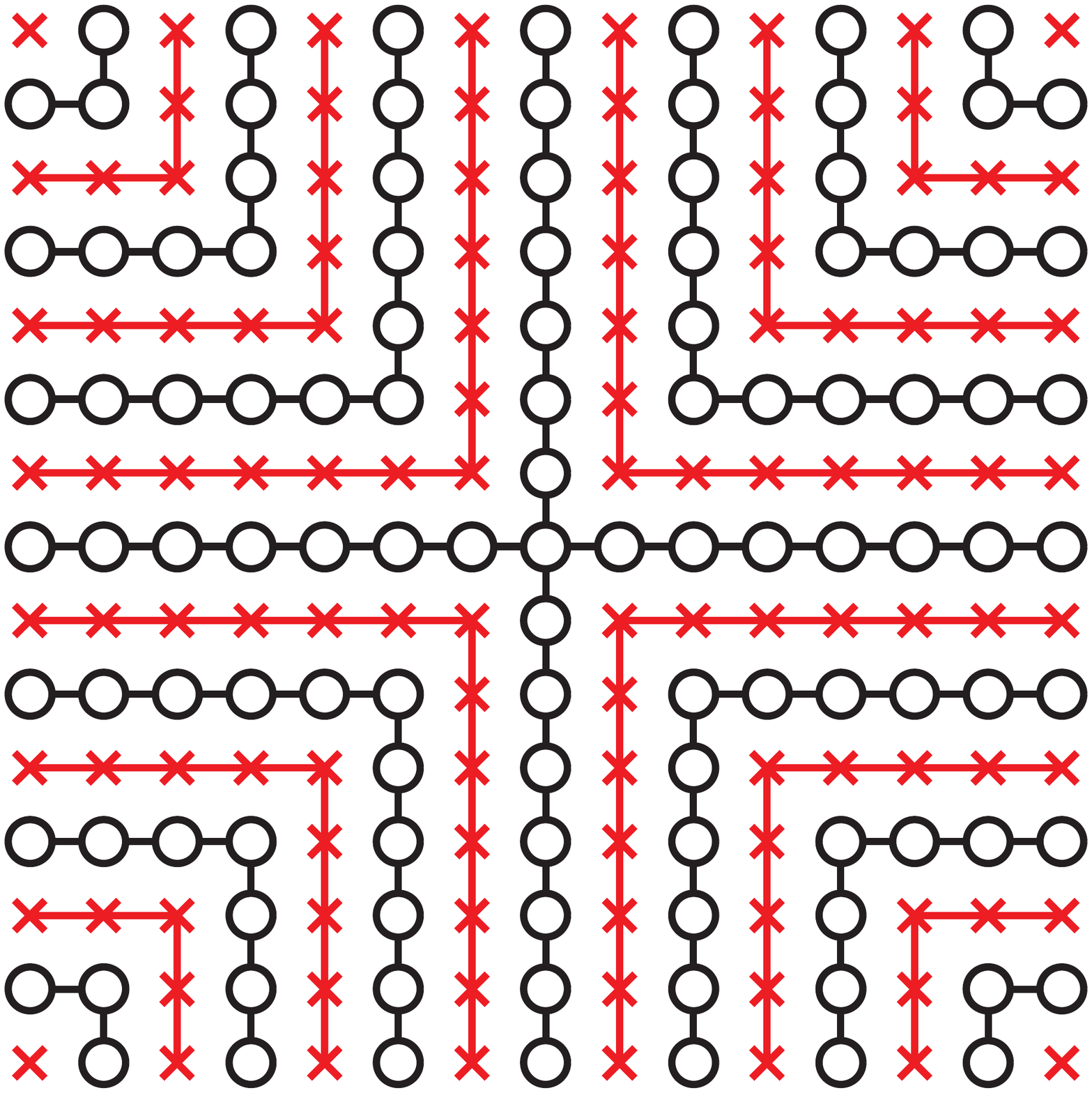}} \quad \subfloat[][\emph{Rectangular grid ($n_x > n_y$)}]{\includegraphics[height=0.35\columnwidth]{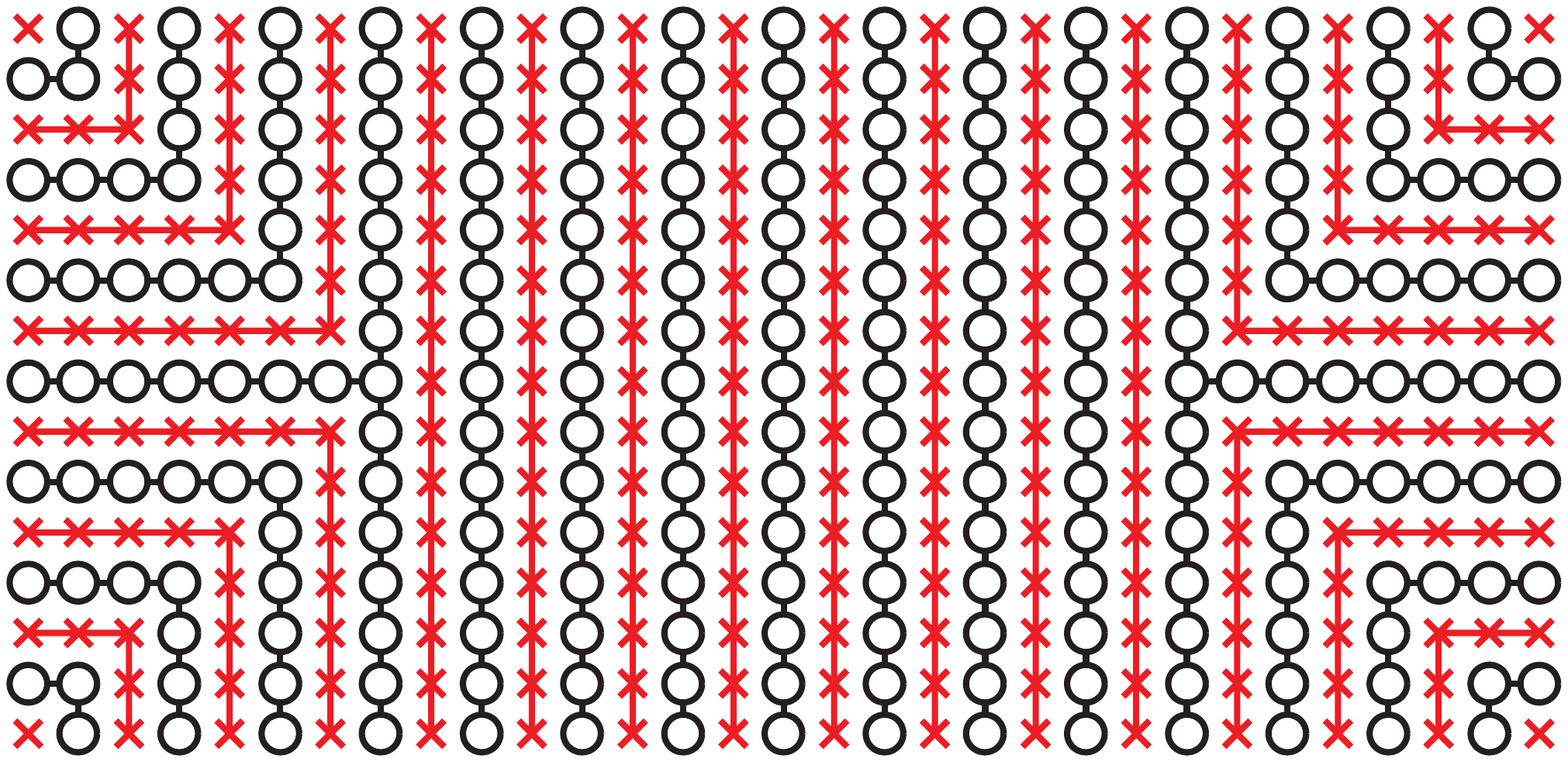}}
	\caption{Relaxation motif for 2D tweed smoothing.}\label{Fig4}
	\subfloat[\emph{Square grid ($n_x = n_y$)}]{\includegraphics[height=0.35\columnwidth]{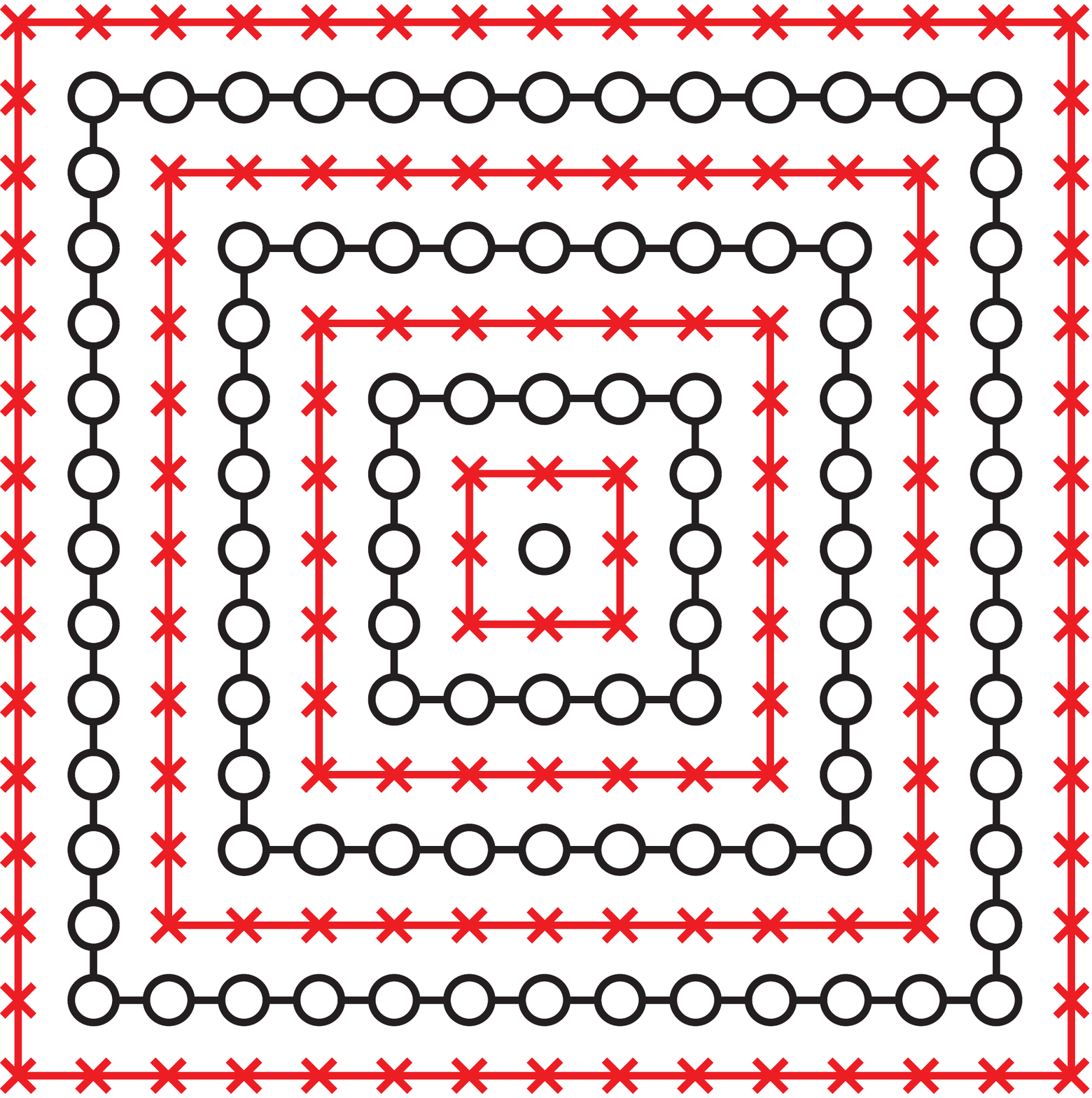}} \quad \subfloat[][\emph{Rectangular grid ($n_x > n_y$)}]{\includegraphics[height=0.35\columnwidth]{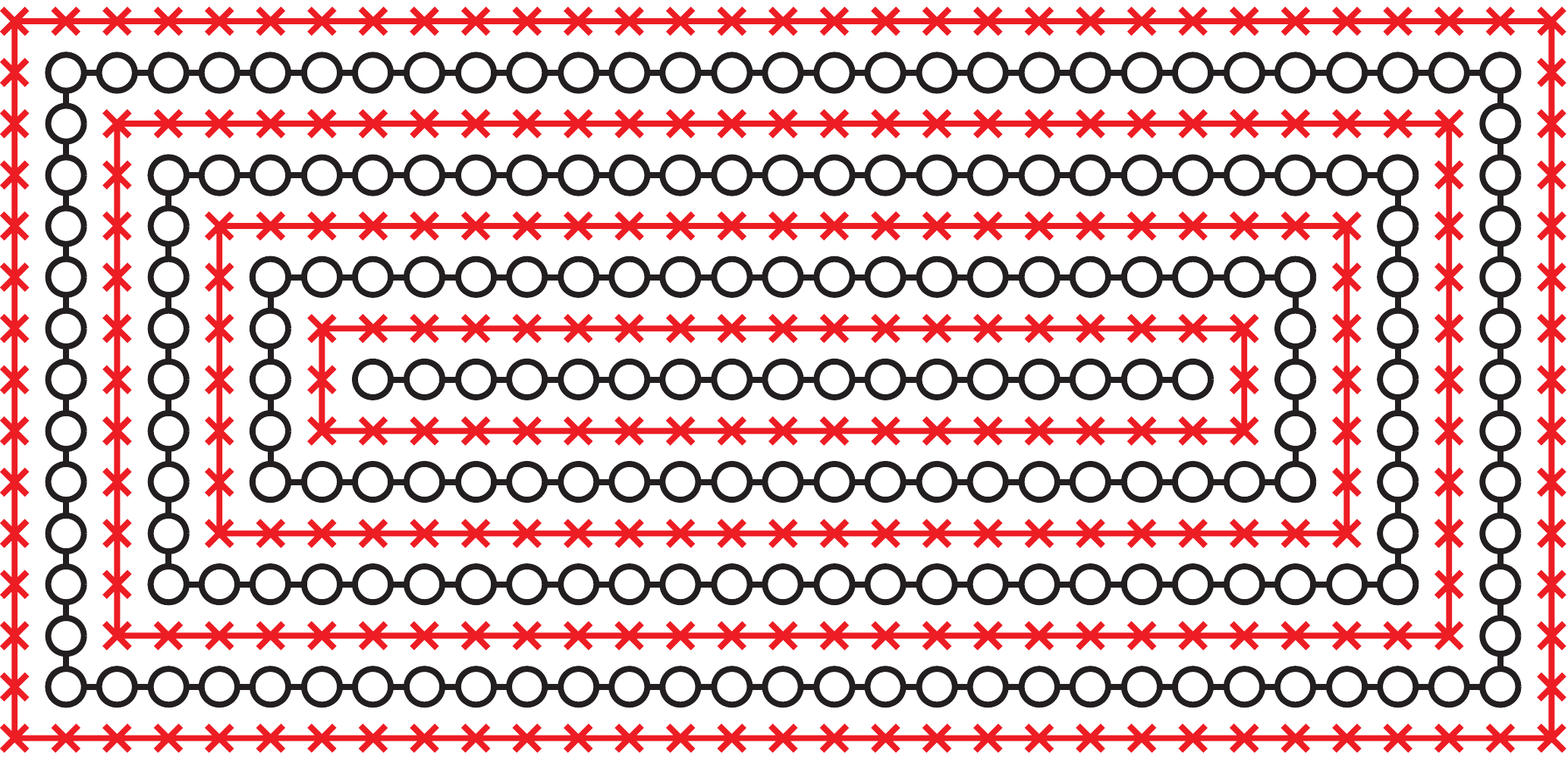}}
	\caption{Relaxation motif for 2D wireframe smoothing.}\label{Fig5}
\end{figure}

Tweed relaxation iteratively solves the linear system of equations arising from the second-order central discretization of second-order 2D or 3D elliptic PDEs.  A prototypical example is given in \eqref{prob:DP}; for convenience, we take $n_x$ and $n_y$ (and, in the 3D case, $n_z$) as even.
As depicted in Figure~\ref{Fig4}, starting from each corner (labelled as ``red''), points are labeled in blocks of alternating colours by forming horizontal and vertical lines of points drawn perpendicular to the domain boundaries,
and extended until such lines connect inside the domain.
Due to the loose visual analogy between such gridpoint arrangements and certain cloth textures, this smoothing scheme has been dubbed ``tweed'' relaxation.
In the 2D case with $n_x = n_y$, as illustrated in Figure~\ref{Fig4}a, four legs of the same colour converge at the centre of the domain.
In the 2D case with $n_x \ne n_y$, as illustrated in Figure~\ref{Fig4}b, two three-legged blocks arise, and the remaining gridpoints in the central part of the domain are connected by lines of alternating colour perpendicular to the closest boundaries.  Again, brief discussion of the natural extension of this motif to the 3D case is deferred to \S \ref{sec:3d}.

Following an approach analogous to that used in zebra relaxation, in tweed relaxation,
\eqref{prob:DP} is first solved exactly at each red point while holding the values of the unknowns at the black points fixed,
then \eqref{prob:DP} is solved exactly at each black point while holding the values of the unknowns at the red points fixed.
In the case of zebra relaxation, such an approach leads, for each (linear) block of a given colour, to a single tridiagonal system of equations that may be solved efficiently using the Thomas algorithm.
In the case of 2D tweed relaxation, this approach leads, for each block other than the corners\footnote{At the corners, simple pointwise relaxation is performed.},
to $m$ tridiagonal systems that are interconnected at a common branch point, where $m=2$, 3, or 4;
an efficient technique to solve this subproblem\footnote{Of course, the case with $m=2$ can be solved directly with the Thomas algorithm itself.}, which we refer to as the {$m$-legged Thomas algorithm}, is described in~\cite{NR}, with a working Matlab implementation provided in ~\cite{NRC}.  In short, the {$m$-legged Thomas algorithm} performs a forward sweep from the tips of each leg in towards the branch point, then performs a solve relating the branch point to the points nearest to the branch point along each leg,
and finally performs a back substitution going back out to the tips of each leg.

\begin{figure}  
	\centering
	\psfrag{x}[t][c]{\small$x$}\psfrag{y}[b][c]{\small$y$}\includegraphics[width=0.45\columnwidth]{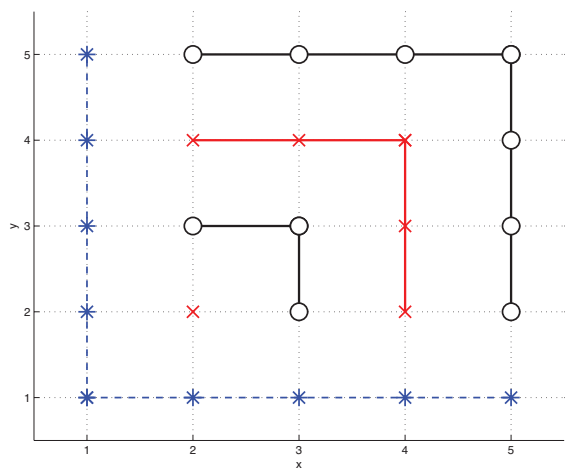}
	\caption{Tweed motif near a corner of a 2D domain. Blue dashed lines indicate the domain boundaries where the value of the unknown is specified.}
	\label{Fig6}
\end{figure}

To illustrate, consider the iterative solution of \eqref{prob:DP} over a uniform grid; defining $C=-2/\Delta x^2-2/\Delta y^2$, $EW=1/\Delta x^2$, and $NS=1/\Delta y^2$,
the block of red points that connect at the $(4,\,4)$ gridpoint in Figure~\ref{Fig6} in this case are governed by the following equations:
\begin{equation*}
\small
\begin{aligned}
	\begin{bmatrix}
		1 & & & \\
		EW & C & EW & \\
		& EW & C & EW
	\end{bmatrix} \,
	\begin{pmatrix}
		u_{1,\,4} \\
		u_{2,\,4} \\
		u_{3,\,4}
	\end{pmatrix} &= 
	\begin{pmatrix}
		g_{1,\,4} \\
		f_{2,\,4} - (u_{2,\,3} + u_{2,\,5}) / \Delta y^2 \\
		f_{3,\,4} - (u_{3,\,3} + u_{3,\,5}) / \Delta y^2
	\end{pmatrix}, \\[3pt]
	\begin{bmatrix}
		1 & & & \\
		NS & C & NS & \\
		& NS & C & NS
	\end{bmatrix}\,
	\begin{pmatrix}
		u_{4,\,1} \\
		u_{4,\,2} \\
		u_{4,\,3}
	\end{pmatrix} &= 
	\begin{pmatrix}
		g_{4,\,1} \\
		f_{4,\,2} - (u_{3,\,2} + u_{5,\,2}) / \Delta x^2 \\
		f_{4,\,3} - (u_{3,\,3} + u_{5,\,3}) / \Delta x^2
	\end{pmatrix}, \\[3pt]
	\begin{bmatrix}
		EW & NS & C 
	\end{bmatrix}\,
	\begin{pmatrix}
	u_{3,\,4} \\
	u_{4,\,3} \\
	u_{4,\,4}
	\end{pmatrix} &=
	f_{4,\,4} - u_{5,\,4} / \Delta x^2 + u_{4,\,5} / \Delta y^2.
\end{aligned}
\end{equation*}
More generally, each tweed relaxation involving $m$ legs of length $p+1$ (including the branch point) can be formulated as the solution of a system of linear equations of the following form:
\begin{subequations}\label{eq:mLThomas}
\begin{align}
	\begin{bmatrix}
		b_1^{(i)} & c_1^{(i)} & & & & \\
		a_2^{(i)} & b_2^{(i)} & c_2^{(i)} & & & \\
		& \ddots & \ddots & \ddots & & \\
		& & a_{p-1}^{(i)} & b_{p-1}^{(i)} & c_{p-1}^{(i)} & \\
		& & & a_p^{(i)} & b_p^{(i)} & c_p^{(i)}		
	\end{bmatrix} \,
	\begin{pmatrix}
		x_1^{(i)} \\
		x_2^{(i)} \\
		\vdots \\
		x_{p-1}^{(i)} \\
		x_p^{(i)} \\
		x_{\textrm{centre}}
	\end{pmatrix} &=
	\begin{pmatrix}
		r_1^{(i)} \\
		r_2^{(i)} \\
		\vdots \\
		r_{p-1}^{(i)} \\
		r_p^{(i)}
	\end{pmatrix} \quad i = 1,\,2,\,\dots, m, \label{eq:mLThomas;a} \\[3pt]
	\begin{bmatrix}
		d_1 & d_2 & \dots & d_m & d_{m+1}
	\end{bmatrix} \,
	\begin{pmatrix}
		x_p^{(1)} \\
		x_p^{(2)} \\
		\vdots \\
		x_p^{(m)} \\
		x_{\textrm{centre}}
	\end{pmatrix} &= r_{\textrm{centre}}; \label{eq:mLThomas;b}
\end{align}
\end{subequations}
note that, on a stretched grid, the computations of the RHS terms $r^i_j$ along the legs typically require 4 floating-point operations (flops).
This system of $mp+1$ equations in $mp+1$ unknowns can be solved efficiently by applying the \emph{$m$-legged} extension of the Thomas algorithm described previously. The system of equations in~\eqref{eq:mLThomas} requires $8mp+m+1$ flops to solve.

Tweed relaxation applied to the iterative solution of \eqref{prob:DP} over an $n \times n$ square grid requires, at each relaxation iteration, $4$ pointwise relaxations (one for each corner), $2(n-3)$ two-legged Thomas solves with legs of increasing size starting from the corners of the domain towards the centre, and one four-legged Thomas relaxation for the centre cross. Hence, the number of flops required for applying a full round of tweed relaxations is:
\begin{itemize}\itemsep-2pt
	\item Corners: $4 \times 9 = 36$ \textit{(Pointwise relaxation)}
	\item Corner legs: $4 \times \sum\limits_{i=1}^{(n-3)/2} \left[ (8\,i+ 4) \, \text{\textit{(RHS computation)}} + \textsc{mLThomas}(p = i,\,m = 2) \right] = 12\,n^2 - 34\,n - 6$
	\item Centre cross: $4 \times 2(n-1) \, \text{\textit{(RHS computation)}} + 1 \times \textsc{mLThomas}(p=(n-1)/2,\,m=4) = 24\,n - 19$
\end{itemize}
Overall, $12\,n^2 - 10\,n + 11$ flops are needed (to leading order, taking $N=n^2$, $\sim 8N$ for the forward sweeps and back substitutions, as in the regular Thomas algorithm, and $\sim 4N$ for the RHS computations).

The 2D wireframe relaxation strategy starts by performing a block relaxation on all points adjacent to the boundaries, and proceeds towards the centre while alternating the colour of the blocks. This creates a pattern of concentric box-shaped blocks of alternate colour, as depicted in Figure~\ref{Fig5}. The relaxation scheme terminates with one pointwise relaxation if $n_x = n_y$, as shown in Figure~\ref{Fig5}a, or with one line relaxation in the $x$ direction if $n_x > n_y$, as shown in Figure~\ref{Fig5}b.

\begin{figure}  
	\centering
	\psfrag{x}[t][c]{\small$x$}\psfrag{y}[b][c]{\small$y$}\includegraphics[width=0.45\columnwidth]{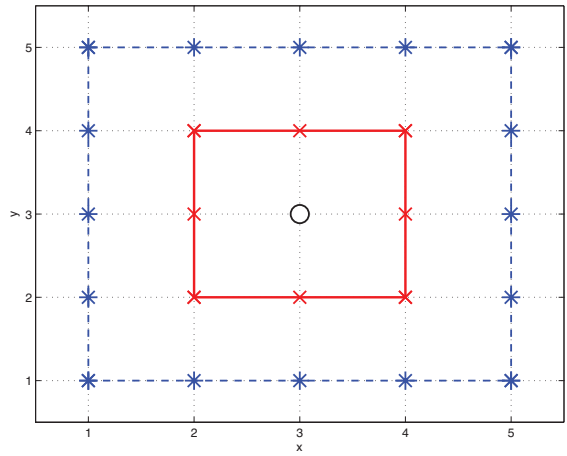}
	\caption{Wireframe relaxation scheme for a $5 \times 5$ uniform grid on a Cartesian domain. Blue dashed lines indicate the domain boundaries, red lines connect the gridpoints involved in the same block relaxation.}
	\label{Fig7}
\end{figure}

To illustrate, consider the iterative solution of \eqref{prob:DP} over a uniform grid, with a red box with corners $(2,\,2),\,(2,\,4),\,(4,\,2),\,(4,\,4)$, as depicted in Figure~\ref{Fig7}.
The tridiagonal circulant system of equations associated with such a relaxation is:
\begin{gather*}
	\begin{pmatrix}
		C & EW & & & & & & NS \\
		EW & C & EW & & & & & \\
		& EW & C & NS & & & & \\
		& & NS & C & NS & & & \\
		& & & NS & C & EW & & \\
		& & & & EW & C & EW & \\
		& & & & & EW & C & NS \\
		NS & & & & & & NS & C
	\end{pmatrix}
	\begin{pmatrix}
		u_{2,\,2} \\
		u_{3,\,2} \\
		u_{4,\,2} \\
		u_{4,\,3} \\
		u_{4,\,4} \\
		u_{3,\,4} \\
		u_{2,\,4} \\
		u_{2,\,3}
	\end{pmatrix} = 
	\begin{pmatrix}
		f_{2,\,2} - u_{1,\,2} / \Delta x^2 - u_{2,\,1} / \Delta y^2 \\
		f_{3,\,2} - u_{3,\,1} / \Delta y^2 - u_{3,\,3} / \Delta y^2 \\
		f_{4,\,2} - u_{4,\,1} / \Delta y^2 - u_{5,\,2} / \Delta x^2 \\
		f_{4,\,3} - u_{3,\,3} / \Delta x^2 - u_{5,\,3} / \Delta x^2 \\
		f_{4,\,4} - u_{5,\,4} / \Delta x^2 - u_{4,\,5} / \Delta y^2 \\
		f_{3,\,4} - u_{3,\,5} / \Delta y^2 - u_{3,\,3} / \Delta y^2 \\
		f_{2,\,4} - u_{2,\,5} / \Delta y^2 - u_{1,\,4} / \Delta x^2 \\
		f_{2,\,3} - u_{1,\,3} / \Delta x^2 - u_{3,\,3} / \Delta x^2
	\end{pmatrix}
\end{gather*}
More generally, each wireframe relaxation involving a block of $n$ points connected together requires the solution of the following linear system:
\begin{equation}
	\begin{bmatrix}
		b_1  & c_1   &              &             & a_1       \\
		a_2  & b_2   & c_2       &             &              \\
		        & \ddots & \ddots   & \ddots &              \\
		        &          & a_{m-1} & b_{m-1} & c_{m-1} \\
		 c_m &          &              & a_m      & b_m 
	\end{bmatrix} \,
	\begin{pmatrix}
		x_1 \\
		x_2 \\
		\vdots \\
		x_{m-1} \\
		x_m
	\end{pmatrix} = 
	\begin{pmatrix}
		r_1 \\
		r_2 \\
		\vdots \\
		r_{m-1} \\
		r_m
	\end{pmatrix} \label{eq:CircTri}
\end{equation}
The circulant tridiagonal system in~\eqref{eq:CircTri} can be solved using a periodic (a.k.a.~circulant) Thomas solver, as introduced in~\cite{ahlberg1967}. A minimal storage implementation of such algorithm is presented in~\cite{cavaglieri2016new}. For an $m \times m$ system like that in~\eqref{eq:CircTri}, the circulant Thomas solver requires $14m-16$ flops, which is about $14/8=1.75$ times the computational cost of the Thomas algorithm for tridiagonal systems.

Wireframe relaxation applied to the iterative solution of \eqref{prob:DP} over an $n \times n$ square grid requires, at each iteration, $1$ pointwise relaxation at the centre point and $(n-1)/2$ concentric box relaxations. Hence, the number of flops required for applying a full round of wireframe relaxation is:
\begin{itemize}\itemsep-2pt
	\item Centre point: 9 \textit{(Pointwise relaxation)}
	\item Elsewhere: $\sum\limits_{i = 1}^{(n-1)/2} \left[ 32\,i \, \text{\textit{(RHS computation)}} + \textsc{CircThomas}(m = 8\,i) \right] = 18\,n^2 - 8\,n - 10$
\end{itemize}
Overall, $18\,n^2 - 8\,n - 1$ flops are needed
(to leading order, taking $N=n^2$, $\sim 14N$ for the circulant Thomas solves, and $\sim 4N$ for the RHS computations), which makes wireframe relaxation about $50\%$ more expensive than tweed relaxation.

As some related points of comparison, the computational cost of other leading 2D relaxation schemes, applied to \eqref{prob:DP} on an $n \times n$ grid, are now reviewed:\\
$\bullet$ For checkerboard relaxation (which works well only on unstretched grids), $\sim 9\,n^2$ flops are required per relaxation iteration, as $9$ flops are required at each gridpoint.\\
$\bullet$ For one-direction zebra relaxation (which works well on grids that are stretched in one spatial direction only), $\sim 12\,n^2$ flops are required per relaxation iteration, as the Thomas algorithm requires $\sim 8\,n$ flops and the computation of the RHS requires $4\,n$ flops on each of $n$ lines.\\
$\bullet$ For alternating-direction zebra relaxation (which works well on grids that are streched in two spatial directions), $\sim 24\,n^2$ flops are required per {\it pair} of relaxation iterations; that is,
$\sim 12\,n^2$ flops are required for each relaxation iteration, and these zebra relaxations must alternate between each coordinate direction.

In short, tweed relaxation has a leading-order computational cost that is the same as one-direction zebra relaxation.  However, to successfully apply zebra relaxation to 2D grids that are stretched in both spatial directions, a {\it pair} of relaxation iterations are required (one in each coordinate direction) at each step in the multigrid algorithm, rendering tweed much more efficient as, we will demonstrate.  

Another approach worth mentioning here is the sequential application of tweed and wireframe relaxation, akin to alternating-direction zebra relaxation approach mentioned above. In the 2D case, alternating tweed and wireframe relaxation iterations requires $\sim 30\,n^2$ flops per pair of relaxation iterations (or, averaged over several steps, $\sim 15\,n^2$ flops per relaxation iteration), which makes it $20\%$ more expensive than tweed relaxation or one direction of zebra relaxation.  Alternating tweed/wireframe relaxation might prove useful in the future, when grid clustering is not isolated to localized regions either near the centre or close to the boundaries, but instead is applied in multiple regions, in different directions, within a computational domain.

The geometric multigrid algorithm (described in \S \ref{sec:Multigrid}) is applied (in \S \ref{sec:spectral} and \S \ref{sec:Results}) to the spatially discretized elliptic PDE \eqref{prob:DP} on a $128 \times 128$ grid with different levels of grid stretching, using the tweed, wireframe, and alternating tweed/wireframe relaxation schemes discussed above.  Performance is compared (both analytically and numerically) with multigrid leveraging the traditional checkerboard, one-direction zebra, and alternating-direction zebra relaxation schemes.  
\section{Summary of the geometric multigrid algorithm}
\label{sec:Multigrid}

Consider the problem of solving \eqref{prob:DP}, of the form
$\mathcal{L} u = f$,
on a stretched $(n_x+1)\times (n_y+1)$ grid $\Omega_p$ (including boundary points) via a multigrid algorithm,
with $n_x=2^{\,p}\, a$ and $n_y=2^{\,p}\, b$ where $a$ and $b$ are small positive integers, at least one of which is odd.  The multigrid algorithm leverages a sequence of grids
$\Omega_p,\Omega_{p-1},\dots,\Omega_0$, where $\Omega_{\ell-1}$ is obtained by coarsening $\Omega_\ell$ by a factor of two in each direction (that is, removing every other interior grid line in each direction), and thus $\Omega_0$ is
$(a+1)\times (b+1)$. We indicate with $\mathcal{L}_\ell$ the discretized Laplacian on $\Omega_\ell$, as defined in \eqref{prob:DP}, and will iterate on a sequence of discretized Poisson problems of the form
$\mathcal{L}_\ell u^\ell = f^{\,\ell}$
on $\Omega_\ell$ for $\ell \in [0,p]$.
A skeleton V-cycle linear multigrid algorithm is composed of the following steps [the restriction step in (3) and the prolongation step in (5) are described further in the paragraphs that follow]:
\begin{enumerate}[(1)]\itemsep-2pt
    \item Initialize $u^p=0$ and $\ell=p$.
    \item Apply $\nu_1$ {\it pre-smoothing} relaxations to the problem $\mathcal{L}_\ell u^\ell = f^{\ell}$ on the grid $\Omega_\ell$.
	\item Compute the remaining {\it defect} $d^\ell = f^\ell - \mathcal{L}_\ell u^\ell$ of the solution $u^\ell$, and {\it restrict} this defect $d^\ell$ from $\Omega_\ell$ to the next coarser grid $\Omega_{\ell-1}$, calling the result $f^{\ell-1}$. Set $\ell\leftarrow \ell-1$.
	\item If $\ell>0$, initialize $u^\ell=0$ and repeat from (2); otherwise, solve
	$\mathcal{L}_0 u^0 = f^0$ for the correction $u^0$ directly.
	\item {\it Prolongate} the correction $u^\ell$ on $\Omega_{\ell}$ to the next finer grid $\Omega_{\ell+1}$, calling the result $v^{\ell+1}$, and {\it update} the solution $u^{\ell+1}$ defined on $\Omega_{\ell+1}$ with the correction $v^{\ell+1}$; that is, take
	$u^{\ell+1} \leftarrow u^{\ell+1} + v^{\ell+1}$.  Set $\ell \leftarrow \ell+1$.
	\item Apply $\nu_2$ {\it post-smoothing} relaxations to the problem $\mathcal{L}_\ell u^\ell = f^{\ell}$ on the grid $\Omega_\ell$.  If $\ell<p$, repeat from (5); otherwise, repeat from (2) until the norm of the defect on $\Omega_p$,
	$||f^p - \mathcal{L}_p u^p ||$, is sufficiently small.
\end{enumerate}
The reason that the multigrid algorithm works so well is that most effective relaxation schemes, such as checkerboard relaxation, {\it smooth} the error in the solution very quickly (that is, they significantly refine the solution vector $u^\ell$ on the smallest scales representable on the grid  $\Omega_\ell$ being used), but are inefficient at reducing the defect on the larger length scales representable on the grid; thus, following the multigrid approach, the larger length scales of the defect are addressed by applying smoothing to successively coarser representations of the problem at hand, as described above.
There are various ways to accelerate the multigrid algorithm.  First, to reduce storage, the computation of the defect and its restriction to the coarser grid, in step (3), can be combined into a single step, thereby eliminating the need for the intermediate storage of $d^\ell$.  Similarly, the computation of the prolongation of the correction, $v^{\ell+1}$, and its use in updating the solution $u^{\ell+1}$, in step (5), can also be combined into a single step, thereby eliminating the need for intermediate storage of $v^{\ell+1}$.  Other cycling strategies, performing more iterations at the coarser length scales, are sometimes used.  

Six different relaxation schemes are considered hereafter for the smoothing applied at steps (2) and (6), namely: checkerboard, one-direction zebra, alternating-direction zebra, tweed, wireframe, and alternating tweed/wireframe. For the restriction step in (3), half weighting and full weighting are considered: denoting with $d^\ell$ the defect on the grid $\Omega_\ell$, and with $f^{\ell-1}$ this defect restricted onto the next coarser grid $\Omega_{\ell-1}$, the half-weighting restriction operation is
\begin{equation*}
	{f}^{\ell-1}_{i,\,j} = \frac{1}{2} d^\ell_{2i,\,2j} + \frac{1}{8} (d^\ell_{2i-1,\,2j} + d^\ell_{2i,\,2j-1} + d^\ell_{2i+1,\,2j} + d^\ell_{2i,\,2j+1}),
\end{equation*}
whereas the full-weighting restriction operation is
\begin{align*}
{f}^{\ell-1}_{i,\,j} = \frac{1}{4} d^\ell_{2i,\,2j} &+ \frac{1}{8} (d^\ell_{2i-1,\,2j} + d^\ell_{2i,\,2j-1} + d^\ell_{2i+1,\,2j} + d^\ell_{2i,\,2j+1}) + \\
&+ \frac{1}{16} (d^\ell_{2i-1,\,2j-1} + d^\ell_{2i+1,\,2j-1} + d^\ell_{2i-1,\,2j+1} + d^\ell_{2i+1,\,2j+1}).
\end{align*}
For the prolongation step in (5), bilinear interpolation is used, which is the dual of the full-weighting restriction operation: denoting with $u^{\ell}$ the correction on the grid $\Omega_\ell$, and with $v^{\ell+1}$ this correction prolongated onto the next finer grid $\Omega_{\ell+1}$, the bilinear interpolation operation is
\begin{equation*}
	v^{\ell+1}_{i,\,j} = \left\{ 
	\begin{array}{lll}
		u^{\ell}_{i/2,\,j/2} & i = \text{even}, & j = \text{even} \\[2pt]
		\frac{1}{2} (u^{\ell}_{i/2,\,(j-1)/2} + u^{\ell}_{i/2,\,(j+1)/2}) & i = \text{even}, & j = \text{odd} \\[2pt]
		\frac{1}{2} (u^{\ell}_{(i-1)/2,\,j/2} + u^{\ell}_{(i+1)/2,\,j/2}) & i = \text{odd}, & j = \text{even} \\[2pt]
		\frac{1}{4} (u^{\ell}_{(i-1)/2,\,(j-1)/2} + u^{\ell}_{(i+1)/2,\,(j-1)/2} + u^{\ell}_{(i-1)/2,\,(j+1)/2} + u^{\ell}_{(i+1)/2,\,(j+1)/2}) & i = \text{odd}, & j = \text{odd}
	\end{array}
	\right.
\end{equation*}

\section{Convergence analysis via calculation of spectral radius}\label{sec:spectral}

Calculation of two-grid convergence factors of the associated multigrid operator provides a useful indication of the effectiveness of the combined application of different restriction and prolongation schemes, smoothers, and the number of pre-smoothing and post-smoothing relaxations applied. To proceed, consider an $(n_x+1) \times (n_y+1)$ grid $\Omega_\ell$ and a coarsened $(n_x/2+1) \times (n_y/2+1)$ grid $\Omega_{\ell-1}$. Four different cases of grid stretching are also considered: a uniform grid with $\Delta x = \Delta y$, and three stretched grids, two exhibiting differing amounts of near-wall clustering, and one exhibiting near-centre clustering (see Figure~\ref{Fig2}). Near-wall clustering of $\Omega_\ell$ over the domain $[0,\,L_x] \times [0,\,L_y]$ is achieved with the hyperbolic tangent stretching function
\begin{equation}
\begin{alignedat}{2}
	x_i &= (L_x / 2) \{1+\tanh{ [ c \,(2\,i/n_x - 1)] }/ \tanh{c}\}, \quad && i=0,\ldots,n_x, \\
	y_j &= (L_y / 2) \{1+\tanh{ [ c \,(2\,j/n_y - 1)] }/ \tanh{c}\}, \quad && j=0,\ldots,n_y,
\end{alignedat}
\label{eq:tanh}
\end{equation}
where $c$ is a tuning parameter that determines the amount of stretching:
$c = 1.5$ creates mild stretching (see Figure~\ref{Fig2}a), and
$c = 3.0$ creates more significant stretching (see Figure~\ref{Fig2}b). Near-centre clustering is achieved with a simple shifted version of the stretching function used in~\eqref{eq:tanh} such that
\begin{equation}
\begin{aligned}
		x_i &=  \begin{cases}
		(L_x/2) \tanh{ [ 2\,c\,i/n_x ] }/ \tanh{c}, \quad & i=0,\ldots,n_x/2, \\
		(L_x/2)\{2-\tanh{ [ c \,(2-2\,i/n_x)] }/  \tanh{c}\}, \quad & i=n_x/2,\ldots,n_x; \end{cases} \\
		y_j &= \begin{cases}
		(L_y/2) \tanh{ [ 2\,c\,j/n_y ] }/ \tanh{c}, \quad & j=0,\ldots,n_y/2, \\
		(L_y/2)\{2-\tanh{ [ c \,(2-2\,j/n_y)] }/  \tanh{c}\}, \quad & j=n_y/2,\ldots,n_y. \end{cases}
\end{aligned}
\label{eq:tanh_inv}
\end{equation}
An example of a stretched grid generated  using~\eqref{eq:tanh_inv}, with $c=1.5$, is shown in Figure~\ref{Fig2}c.

We now denote the restriction operator from the fine to the coarse grid as $I^{\ell-1}_\ell$, and the prolongation operator from the coarse grid to the fine grid as $I^\ell_{\ell-1}$. Pre- and post-smoothing operators on $\Omega_\ell$ are indicated as $S^{\,\nu_1}_\ell$ and $S^{\,\nu_2}_\ell$, where $\nu_1$ and $\nu_2$ indicate the number of times the smoother is applied at each step. Following~\cite{trottenberg2000}, the complete two-grid multigrid operator, denoted $M^{\ell-1}_\ell$, is given by
\begin{equation}
	M^{\ell-1}_\ell = S^{\,\nu_2}_\ell \, (I_\ell - I^\ell_{\ell-1} \mathcal{L}^{-1}_{\ell-1} I^{\ell-1}_\ell \mathcal{L}_\ell) \, S^{\,\nu_1}_\ell, \label{eq:Mh2h}
\end{equation}
where $I_\ell$ is the identity matrix on the fine grid $\Omega_\ell$, and $\mathcal{L}_\ell$ and $\mathcal{L}_{\ell-1}$ denote the discrete Laplace operators on the fine and coarse grids, respectively.

Tables~\ref{tab:Uni} through~\ref{tab:InvStr} show the computation of the spectral radius of the two-grid multigrid operator $M^{\ell-1}_\ell$ applied to the solution of \eqref{prob:DP} over a grid with $n_x = n_y = 128$ and $L_x=L_y=1$ with half-weighting and full-weighting used for the restriction operation. Since only the sum of pre- and post-smoothing steps affects the convergence of the two-grid cycle, the spectral radius is reported as a function of the sum $\nu = \nu_1 + \nu_2$.
\begin{table}[pt]
\caption{Spectral radius of the two-grid multigrid operator $M^{\ell-1}_\ell$ for the homogeneous problem (i.e., a uniform grid) with respect to the number of smoothing steps $\nu$ for different smoothers and restriction operators.}
\label{tab:Uni}
\centering
\begin{tabular}{l|c|c|c|c|c|c|}
\cline{2-7}
& \multicolumn{6}{|c|}{Half-weighting restriction} \\
\cline{2-7}
                      Smoother        & $\nu=1$          & $\nu=2$          & $\nu=3$          & $\nu=4$          & $\nu=5$          & $\nu=6$          \\
\hline
Checkerboard                 & $0.4986$ & $0.1238$ & $0.0340$ & $0.0240$ & $0.0189$ & $0.0155$ \\
\hline
One-direction zebra        & $0.5735$ & $0.5014$ & $0.4935$ & $0.4905$ & $0.4881$ & $0.4858$ \\
\hline
Alternating-direction zebra & $0.5836$ & $0.5912$ & $0.5868$ & $0.5823$ & $0.5778$ & $0.5734$ \\
\hline
Tweed                              & $0.5661$ & $0.5032$ & $0.4951$ & $0.4919$ & $0.4893$ & $0.4869$ \\
\hline
Wireframe  				       & $0.5494$ & $0.5007$ & $0.4929$ & $0.4899$ & $0.4874$ & $0.4850$ \\
\hline
Alternating tweed/wireframe  		       & $0.5833$ & $0.5909$ & $0.5865$ & $0.5820$ & $0.5775$ & $0.5731$ \\
\hline
\end{tabular} \\[15pt]
\begin{tabular}{l|c|c|c|c|c|c|}
\cline{2-7}
& \multicolumn{6}{|c|}{Full-weighting restriction} \\
\cline{2-7}
                   Smoother              & $\nu=1$          & $\nu=2$          & $\nu=3$          & $\nu=4$          & $\nu=5$          & $\nu=6$        \\
\hline
Checkerboard                 &$0.2494$ & $0.0739$ & $0.0526$ & $0.0408$ & $0.0333$ & $0.0283$ \\
\hline
One-direction zebra        & $0.2494$ & $0.0622$ & $0.0167$ & $0.0116$ & $0.0092$ & $0.0077$ \\
\hline
Alternating-direction zebra & $0.0839$ & $0.0391$ & $0.0265$ & $0.0201$ & $0.0162$ & $0.0135$ \\
\hline
Tweed                              & $0.2488$ & $0.0621$ & $0.0163$ & $0.0109$ & $0.0084$ & $0.0069$ \\
\hline
Wireframe				       & $0.2465$ & $0.0612$ & $0.0161$ & $0.0108$ & $0.0084$ & $0.0069$ \\
\hline
Alternating tweed/wireframe  		       & $0.0830$ & $0.0382$ & $0.0257$ & $0.0193$ & $0.0154$ & $0.0128$ \\
\hline
\end{tabular}
\end{table}

\begin{table}[pt]
\caption{Spectral radius of the two-grid multigrid operator $M^{\ell-1}_\ell$ for the inhomogeneous problem with near-wall clustering \eqref{eq:tanh}, taking $c =1.5$, with respect to $\nu$ for different smoothers and restriction operators.}
\label{tab:Str}
\centering
\begin{tabular}{l|c|c|c|c|c|c|}
\cline{2-7}
& \multicolumn{6}{|c|}{Half-weighting restriction} \\
\cline{2-7}
                      Smoother         & $\nu=1$          & $\nu=2$          & $\nu=3$          & $\nu=4$          & $\nu=5$          & $\nu=6$          \\
\hline
Checkerboard                 & $0.7900$ & $0.6244$ & $0.4936$ & $0.3904$ & $0.3090$ & $0.2446$ \\
\hline
One-direction zebra        & $0.7921$ & $0.6302$ & $0.5113$ & $0.4825$ & $0.4774$ & $0.4730$ \\
\hline
Alternating-direction zebra & $0.6879$ & $0.6385$ & $0.6096$ & $0.5903$ & $0.5755$ & $0.5629$ \\
\hline
Tweed                              & $0.5478$ & $0.4969$ & $0.4867$ & $0.4803$ & $0.4744$ & $0.4687$ \\
\hline
Wireframe				       & $0.7921$ & $0.6302$ & $0.5113$ & $0.4840$ & $0.4793$ & $0.4754$ \\
\hline
Alternating tweed/wireframe  		       & $0.6891$ & $0.6490$ & $0.6285$ & $0.6138$ & $0.6013$ & $0.5899$ \\
\hline
\end{tabular} \\[15pt]
\begin{tabular}{l|c|c|c|c|c|c|}
\cline{2-7}
& \multicolumn{6}{|c|}{Full-weighting restriction} \\
\cline{2-7}
                      Smoother         & $\nu=1$          & $\nu=2$          & $\nu=3$          & $\nu=4$          & $\nu=5$          & $\nu=6$ \\
\hline
Checkerboard                  & $0.7855$ & $0.6179$ & $0.4867$ & $0.3840$ & $0.3035$ & $0.2403$ \\
\hline
One-direction zebra         & $0.7863$ & $0.6190$ & $0.4879$ & $0.3851$ & $0.3043$ & $0.2409$ \\
\hline
Alternating-direction zebra & $0.0816$ & $0.0344$ & $0.0218$ & $0.0162$ & $0.0129$ & $0.0107$ \\
\hline
Tweed                              & $0.2108$ & $0.0538$ & $0.0257$ & $0.0184$ & $0.0143$ & $0.0117$ \\
\hline
Wireframe				       & $0.7863$ & $0.6190$ & $0.4879$ & $0.3851$ & $0.3043$ & $0.2409$ \\
\hline
Alternating tweed/wireframe  		       & $0.0696$ & $0.0286$ & $0.0179$ & $0.0126$ & $0.0095$ & $0.0073$ \\
\hline           
\end{tabular}
\end{table}

\begin{table}[pt]
\caption{Spectral radius of the two-grid multigrid operator $M^{\ell-1}_\ell$ for the inhomogeneous problem with near-wall clustering \eqref{eq:tanh}, taking $c=3.0$, with respect to $\nu$ for different smoothers and restriction operators.}
\label{tab:Str2}
\centering
\begin{tabular}{l|c|c|c|c|c|c|}
\cline{2-7}
& \multicolumn{6}{|c|}{Half-weighting restriction} \\
\cline{2-7}
                      Smoother         & $\nu=1$          & $\nu=2$          & $\nu=3$          & $\nu=4$          & $\nu=5$          & $\nu=6$          \\
\hline
Checkerboard                 & $0.9541$ & $0.9103$ & $0.8686$ & $0.8287$ & $0.7907$ & $0.7545$ \\
\hline
One-direction zebra        & $0.9541$ & $0.9104$ & $0.8687$ & $0.8290$ & $0.7911$ & $0.7550$ \\
\hline
Alternating-direction zebra & $0.6356$ & $0.5765$ & $0.5395$ & $0.5092$ & $0.4821$ & $0.4572$ \\
\hline
Tweed                              & $0.5462$ & $0.4799$ & $0.4566$ & $0.4386$ & $0.4220$ & $0.4062$ \\
\hline
Wireframe				       & $0.9541$ & $0.9104$ & $0.8687$ & $0.8290$ & $0.7911$ & $0.7550$ \\
\hline
Alternating tweed/wireframe  		       & $0.6538$ & $0.6087$ & $0.5743$ & $0.5440$ & $0.5161$ & $0.4900$ \\
\hline
\end{tabular} \\[15pt]
\begin{tabular}{l|c|c|c|c|c|c|}
\cline{2-7}
& \multicolumn{6}{|c|}{Full-weighting restriction} \\
\cline{2-7}
                      Smoother         & $\nu=1$          & $\nu=2$          & $\nu=3$          & $\nu=4$          & $\nu=5$          & $\nu=6$ \\
\hline
Checkerboard                  & $0.9534$ & $0.9090$ & $0.8666$ & $0.8263$ & $0.7879$ & $0.7512$ \\
\hline
One-direction zebra         & $0.9534$ & $0.9091$ & $0.8668$ & $0.8265$ & $0.7880$ & $0.7514$ \\
\hline
Alternating-direction zebra & $0.1002$ & $0.0372$ & $0.0219$ & $0.0148$ & $0.0110$ & $0.0088$ \\
\hline
Tweed                              & $0.1866$ & $0.0537$ & $0.0282$ & $0.0204$ & $0.0158$ & $0.0135$ \\
\hline
Wireframe				       & $0.9534$ & $0.9091$ & $0.8668$ & $0.8265$ & $0.7880$ & $0.7514$ \\
\hline
Alternating tweed/wireframe  		       & $0.0987$ & $0.0362$ & $0.0223$ & $0.0161$ & $0.0125$ & $0.0103$ \\
\hline
\end{tabular}
\end{table}

\begin{table}[pt]
\caption{Spectral radius of the two-grid multigrid operator $M^{\ell-1}_\ell$ for the inhomogeneous problem with near-centre clustering \eqref{eq:tanh_inv}, taking $c =1.5$, with respect to $\nu$ for different smoothers and restriction operators.}
\label{tab:InvStr}
\centering
\begin{tabular}{l|c|c|c|c|c|c|}
\cline{2-7}
& \multicolumn{6}{|c|}{Half-weighting restriction} \\
\cline{2-7}
                      Smoother         & $\nu=1$          & $\nu=2$          & $\nu=3$          & $\nu=4$          & $\nu=5$          & $\nu=6$          \\
\hline
Checkerboard                 & $0.8865$ & $0.7859$ & $0.6968$ & $0.6178$ & $0.5478$ & $0.4859$ \\
\hline
One-direction zebra        & $0.8876$ & $0.7892$ & $0.7047$ & $0.6356$ & $0.5822$ & $0.5431$ \\
\hline
Alternating-direction zebra & $0.7918$ & $0.7463$ & $0.7149$ & $0.6911$ & $0.6724$ & $0.6574$ \\
\hline
Tweed                              & $0.8876$ & $0.7892$ & $0.7047$ & $0.6356$ & $0.5823$ & $0.5434$ \\
\hline
Wireframe				       & $0.5496$ & $0.4974$ & $0.4913$ & $0.4888$ & $0.4867$ & $0.4846$ \\
\hline
Alternating tweed/wireframe  		       & $0.5312$ & $0.5433$ & $0.5389$ & $0.5353$ & $0.5320$ & $0.5289$ \\
\hline
\end{tabular} \\[15pt]
\begin{tabular}{l|c|c|c|c|c|c|}
\cline{2-7}
& \multicolumn{6}{|c|}{Full-weighting restriction} \\
\cline{2-7}
                      Smoother         & $\nu=1$          & $\nu=2$          & $\nu=3$          & $\nu=4$          & $\nu=5$          & $\nu=6$         \\
\hline
Checkerboard                  & $0.8826$ & $0.7793$ & $0.6884$ & $0.6084$ & $0.5380$ & $0.4759$ \\
\hline
One-direction zebra         & $0.8829$ & $0.7798$ & $0.6891$ & $0.6092$ & $0.5388$ & $0.4767$ \\
\hline
Alternating-direction zebra & $0.0805$ & $0.0375$ & $0.0249$ & $0.0185$ & $0.0146$ & $0.0119$ \\
\hline
Tweed                              & $0.8829$ & $0.7798$ & $0.6891$ & $0.6092$ & $0.5388$ & $0.4767$ \\
\hline
Wireframe				       & $0.1887$ & $0.0408$ & $0.0257$ & $0.0195$ & $0.0160$ & $0.0136$ \\
\hline
Alternating tweed/wireframe  		       & $0.0805$ & $0.0363$ & $0.0248$ & $0.0188$ & $0.0150$ & $0.0125$ \\
\hline
\end{tabular}
\end{table}

In the homogeneous case (Table~\ref{tab:Uni}),
an appreciably reduced spectral radius
(and, thus, an appreciably improved convergence rate)
is obtained using any of the block relaxation methods, but
only when full-weighting restriction is applied.
[In contrast, using the Checkboard relaxation method,
the half-weighting restriction approach actually works a tad better than the (more-expensive) full-weighting restriction approach.]
Of course, in all cases, convergence improves with the number of smoothing steps $\nu$ applied.

In the stretched grid cases (Tables~\ref{tab:Str} through \ref{tab:InvStr}),
an appreciably reduced spectral radius
(and, thus, an appreciably improved convergence rate)
is obtained using only certain choices of block relaxation methods (discussed further below), but (again) only when
full-weighting restriction is applied.  In these cases, the
Checkboard relaxation method proves to be less effective. 

For the cases with near-wall clustering (Tables~\ref{tab:Str} and~\ref{tab:Str2}), both checkerboard and one-direction zebra exhibit a significant degradation in convergence for both choices of the restriction operator. However, convergence remains fast when alternating-direction zebra or tweed is used for the smoothing, and full-weighting restriction is implemented.
Note in particular that the convergence rate using a pair of alternating-direction zebra relaxations (one in each coordinate direction) is slightly better than that of a single tweed relaxation for small $\nu$, whereas the convergence rate is comparable for higher $\nu$.
Note again, however, that the computational cost of a single tweed relaxation is roughly half the computational cost of a pair of alternating-direction zebra relaxations, thereby rendering tweed relaxation with full-weighting restriction the clearly superior choice.
Further, comparing Tables~\ref{tab:Str} and~\ref{tab:Str2}, it is seen that the convergence of tweed relaxation is affected only slightly by the degree of grid stretching applied.

For the case with near-centre clustering (Table~\ref{tab:InvStr}), again, both checkerboard and one-direction zebra prove to be inadequate.  In this case, convergence remains fast when alternating-direction zebra or wireframe is used for the smoothing, and full-weighting restriction is implemented.
Note in particular that the convergence rate using a pair of alternating-direction zebra relaxations (one in each coordinate direction) is only slightly better than that of a single wireframe relaxation for all values of $\nu$.
Note again, however, that the computational cost of a single wireframe relaxation is about 25\% less than the computational cost of a pair of alternating-direction zebra relaxations, thereby rendering wireframe relaxation with full-weighting restriction a very competitive choice.

Of course, tweed relaxation alone is poorly suited for the near-centre clustering case, and wireframe relaxation alone is poorly suited for the near-wall clustering case.  Alternating tweed and wireframe relaxation steps provides similar performance, per pair of relaxations as tweed alone (in the near-wall clustering case) or wireframe alone (in the near-centre clustering case), but at increased computational cost per relaxation.

Further insight on checkerboard and zebra relaxation may be achieved by rigorous or local Fourier analysis (see \cite{trottenberg2000}). However, the complicated arrangement of gridpoints in tweed and wireframe relaxation prevents the extension of these analysis tools to the new smoothers proposed here.

\section{Overall multigrid performance using tweed and wireframe relaxation}\label{sec:Results}
To assess the overall performance of the 2D tweed and wireframe relaxation schemes introduced in the present work, we applied the multigrid algorithm described in Section~\ref{sec:Multigrid} to the solution of \eqref{prob:DP} over uniform and stretched grids with the different smoothing schemes implemented. The RHS vector $f^\ell$ used in these tests is defined using uniformly-generated random numbers, and full-weighting restriction is used in every simulation reported.  We selected $L_x=L_y=1$, $n_x=n_y=128$ (that is, $p=7$), and $\nu_1 = \nu_2$ (that is, the same number of pre-smoothing and post-smoothing relaxations are used) for every simulation reported.

\begin{figure}  
\centering
\includegraphics[trim=3 0 40 0,clip,width=0.33\columnwidth]{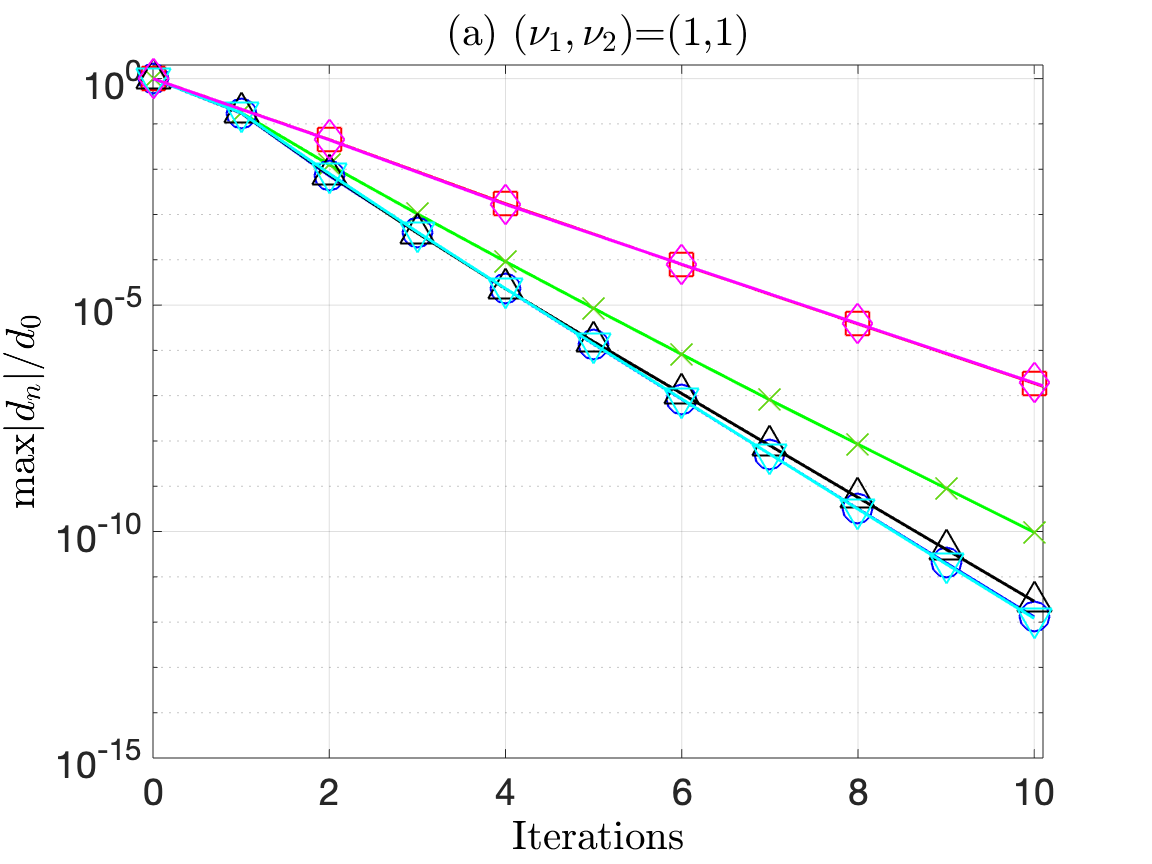}\includegraphics[trim=3 0 40 0,clip,width=0.33\columnwidth]{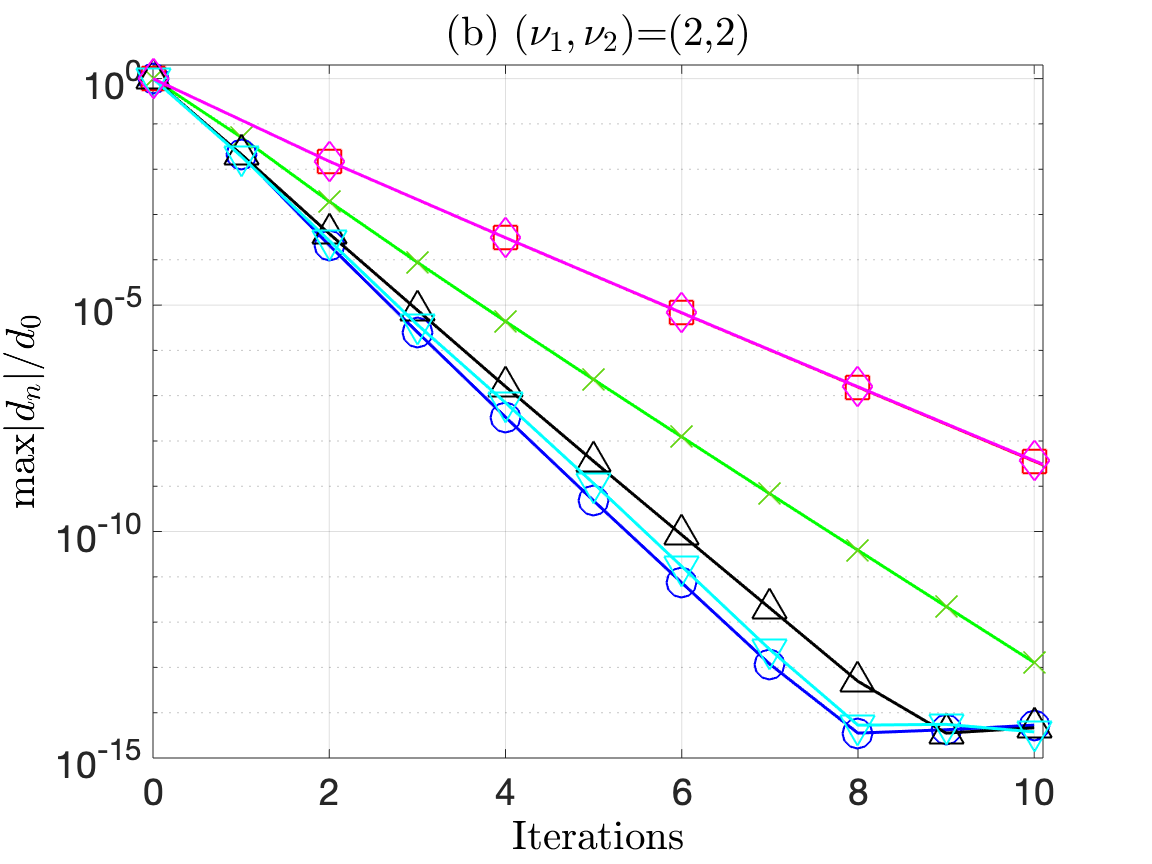}\includegraphics[trim=3 0 40 0,clip,width=0.33\columnwidth]{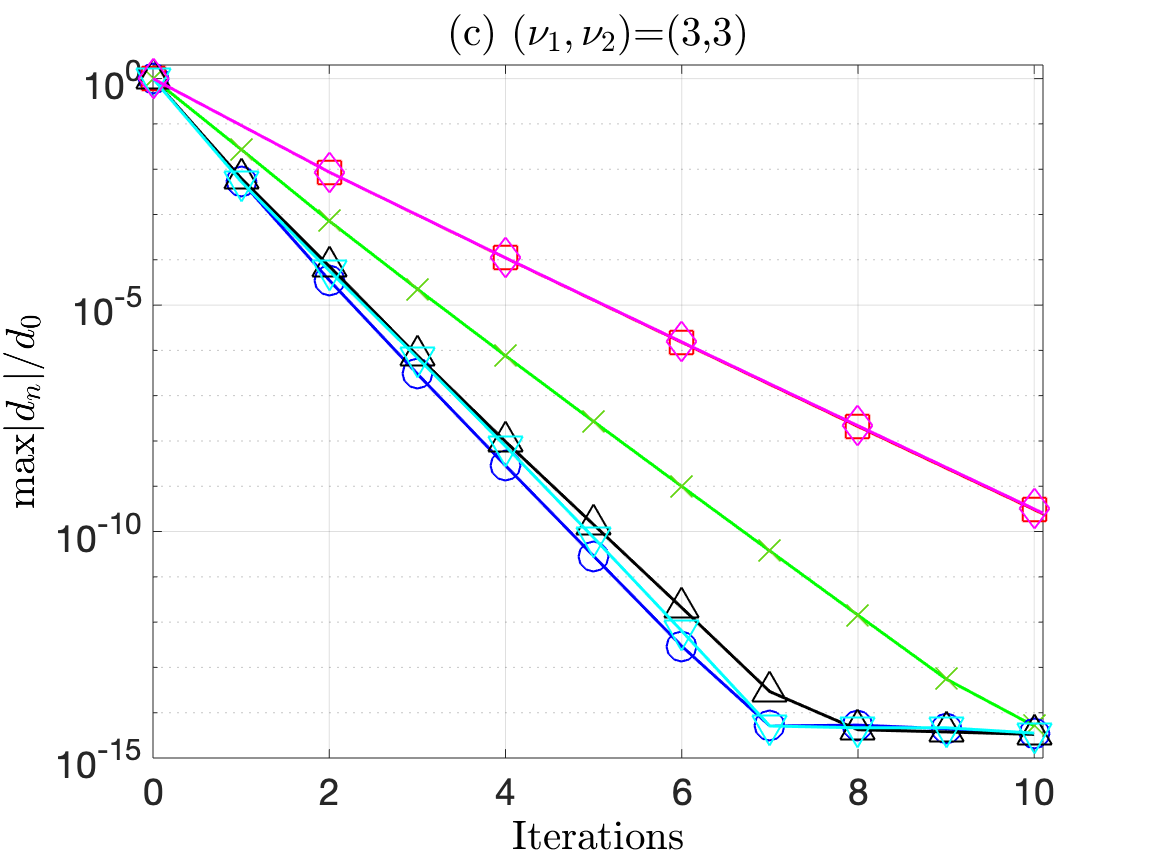}
\caption{Multigrid convergence on \eqref{prob:DP} over a $129 \times 129$ uniform grid with different smoothers: checkerboard (green crosses), one-direction zebra (blue circles), alternating-direction zebra (red squares), tweed (black upward-pointing triangles), wireframe (light blue downward-pointing triangles), alternating tweed/wireframe (magenta diamonds).  To represents these convergence results equitably, all curves are plotted plotted per relaxation on the horizontal axis; note that alternating zebra and alternating tweed/wireframe actually perform two relaxations per iteration, one of each type.}
\label{Fig8}
\end{figure}

Convergence of the multigrid algorithm in the uniform-grid case is reported in Figure~\ref{Fig8}, where the maximum defect, normalized by the initial maximum defect $d_0$, is reported at each multigrid iteratiomn. It can be observed that checkerboard smoothing provides very rapid convergence; though some gains can be seen by introducing various block relaxation schemes, these gains are more than offset (in the uniform-grid case) by the significantly increased computational cost per relaxation associated with these block relaxation schemes. 

Convergence of the multigrid algorithm in the near-wall clustering case, with the grid generated using \eqref{eq:tanh} for $c = 1.5$ and $c=3.0$, is reported in Figure~\ref{Fig9}.
Generally speaking, simulations show good agreement with the theoretical results presented in Section~\ref{sec:spectral}.
In particular, convergence of checkerboard, one-direction zebra, and wireframe relaxation are significantly degraded in the case of $c=1.5$, and actually fail to converge in the case of $c=3$.
The convergence rates (per relaxation step) obtained in the tests performed over a uniform grid (Figure~\ref{Fig8}) are recovered in the near-wall clustering case when alternating-direction zebra is applied, a result that is well-known in the multigrid literature~\cite{trottenberg2000}.
Remarkably, tweed relaxation achieves a level of convergence in five relaxation steps (again, each of which costing $\sim 12\,n^2$ flops) what alternating-direction zebra achieves in five {\it pairs} (that is, ten) relaxation steps (again, each one of which costing $\sim 12\,n^2$ flops).
We thus see that tweed significantly outperforms all other competing block smoothers in the near-wall clustering case. 

As shown in Figure~\ref{fig:defect} for the near-wall clustering case, after three multigrid cycles, the distribution of the absolute value of the defect $d^p$ on $\Omega^p$ is focused near the centre of the domain (i.e., in the region where the grid is coarsest) in the case of tweed, and is distributed more uniformly throughout the domain in the case of alternating-direction zebra.

Convergence of the multigrid algorithm in the near-centre clustering case, with the grid generated using \eqref{eq:tanh_inv} for $c = 1.5$, is reported in Figure~\ref{Fig10}.
Again, simulations show close agreement with the theoretical results presented in Section~\ref{sec:spectral}.
In particular, convergence of checkerboard, one-direction zebra, and tweed relaxation are significantly degraded.
The convergence rates (per relaxation step) are improved when alternating-direction zebra is applied, albeit not recovering the convergence rate in the uniform grid case.
Wireframe relaxation is seen to outperforms  alternating-direction zebra for all cases reported.

As shown in Figure~\ref{fig:defectINV} for the near-centre clustering case, after three multigrid cycles, the distribution of the absolute value of the defect $d^p$ on $\Omega^p$ is focused away from the centre of the domain (where the grid is densest) in both the wireframe and alternating-direction zebra cases. It is also noted that the performance of alternating-direction zebra, wireframe, and alternating tweed/wireframe degrade as $c$ is increased, failing to achieve convergence in the near-centre clustering case when extreme stretching is applied (that is, for  $c \gtrsim 2.3$; not shown here).

\section{Extension of tweed and wireframe relaxation to 3D}\label{sec:3d}

The extension of 2D tweed relaxation to 3D is straightforward, as illustrated in Figure~\ref{Fig13}a. As in the 2D case, the (now, eight) corner points are first relaxed using a pointwise smoother. Then, starting from the corners, points are relaxed in blocks of alternate colours, each composed of $m=2$ or 3 legs.
If $n_x=n_y=n_z$, as illustrated in Figure~\ref{Fig13}a, a single $m=6$ block arises in the centre.
In all cases, it is seen that each gridpoint is a member of exactly one relaxation block.  Further, for large grids, the number of points on the legs dominates the number of branch points.  Thus, to leading order, the computation cost of 3D tweed for large grids is the same as that of a {\it single} set of sweeps (that is, in a single direction) of 1D zebra relaxations.  However, on a grid that is stretched in three directions, the 3D alternating-direction zebra scheme requires {\it three} successive sweeps of 1D zebra relaxations, one in each direction. Hence, the leading-order computational cost of 3D tweed relaxation is, to leading order, {\it one third}
that of the three sweeps of the 3D alternating-direction zebra scheme for such problems.

The extension of 2D wireframe relaxation to 3D is illustrated in Figure \ref{Fig13}b. Starting from the faces of the 3D domain, gridpoints are relaxed in blocks covering the faces of concentric 3D rectangular cuboids filling the domain.  Along all the edges of each such rectangular cuboid, a 3D ``wireframe'' extension of the Thomas algorithm is used, as discussed in \cite{cavaglieri2016new}.  Within each of the six faces of each such rectangular cuboid, a sequence of concentric 2D box relaxations is performed, using the circulant Thomas algorithm, as illustrated in Figure \ref{Fig5}.
If $n_x = n_y = n_z$, as illustrated in Figure \ref{Fig13}b, this sequence of relaxations includes a point relaxation at the centre of each face of each cuboid.

\section{Complementarity}\label{sec:Complementarity}

In 2D (see Figures \ref{Fig4}a/\ref{Fig5}a, or Figures \ref{Fig4}b/\ref{Fig5}b), the tweed and wireframe motifs are disjoint (the branched lines on the grid over which block relaxation is performed are nonoverlapping) and {\it complementary}.  That is, a 2D Cartesian grid has 4 lines emanating from every node (towards its 4 nearest neighbors), and the tweed and wireframe motifs, when superposed, incorporate every single wire of the 2D Cartesian grid.   In 3D (see Figures \ref{Fig13}a/b), the tweed and wireframe motifs are also disjoint.  Further, a 3D Cartesian grid has 6 lines emanating from every node (towards its 6 nearest neighbors).  The tweed and wireframe motifs, together with a third motif which we call “complementary wireframe”, then together form a complementary set of 3 motifs, which when superposed incorporate every single wire of the 3D Cartesian grid in a nonoverlapping fashion.  Complementary wireframe connections again cover the surface of concentric cuboids, and are everywhere orthogonal to to the wireframe motif.  This results in a relaxation scheme with the 8 corner points taken as point relaxations, and non-branched lines (over the 3 adjoining faces) that effectively (in a Cartesian fashion) snake around each corner point of each concentric cuboid (rather than lines that snake around each face-centered point of each concentric cuboid, as in the case of the wireframe motif).  In the case of $n_x=n_y=n_z$, the complementary wireframe motif builds out to a set of 6 adjoined “+” shapes centered on each face of each cuboid.

It is not entirely clear what natural grid stretching patterns the complementary wireframe motif would, taken alone, be especially useful for.  However, for ensuring rapid convergence for more general grid stretching patterns, the analog of the triplets of sweeps of 3D alternating zebra relaxation (that is, zebra in each coordinate direction), leveraging the tweed/wireframe idea, is triplets of sweeps, one each of tweed, wireframe, and complementary wireframe. 

\section{Conclusions}
\label{sec:Conclusions}
Two new relaxation schemes appropriate for the smoothing step in geometric multigrid algorithms applied on 2D stretched grids have been introduced. The implementation of multigrid leveraging such smoothers facilitates the efficient solution of large linear (and, ultimately, nonlinear) systems arising from the discretization of elliptic PDEs on grids that are stretched in multiple spatial directions.

\begin{figure}  
\centering
\includegraphics[trim=3 0 40 0,clip,width=0.33\columnwidth]{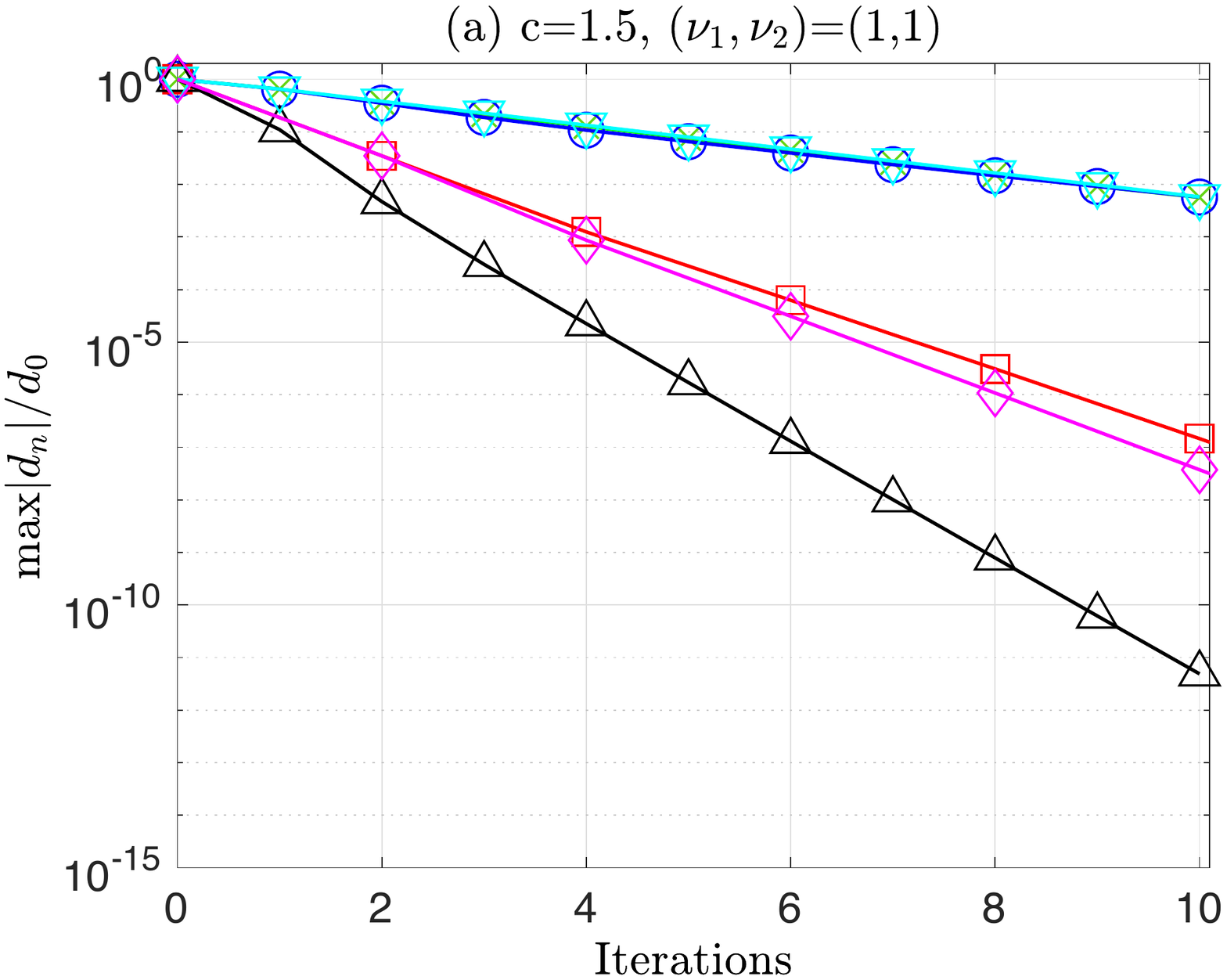}\includegraphics[trim=3 0 40 0,clip,width=0.33\columnwidth]{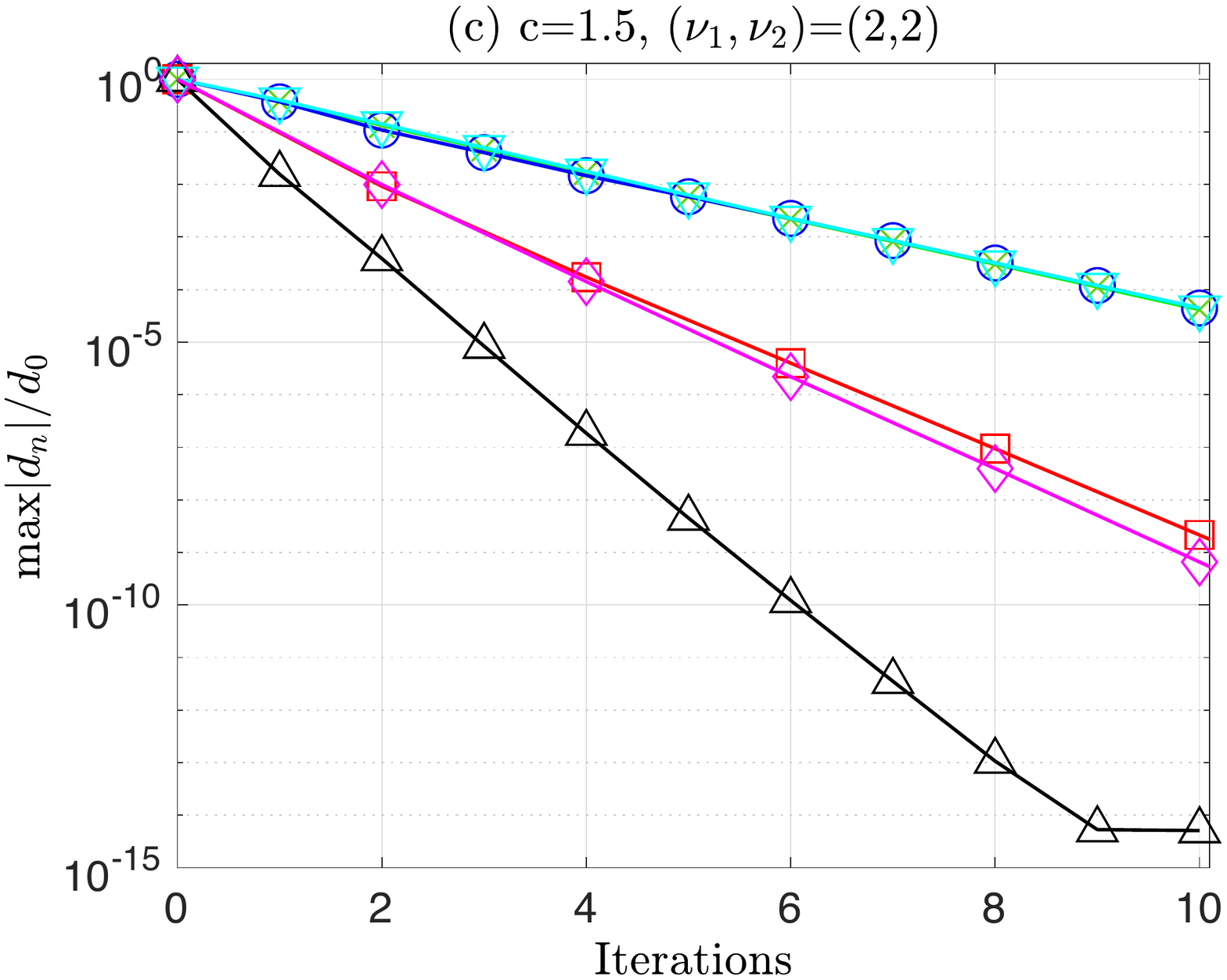}\includegraphics[trim=3 0 40 0,clip,width=0.33\columnwidth]{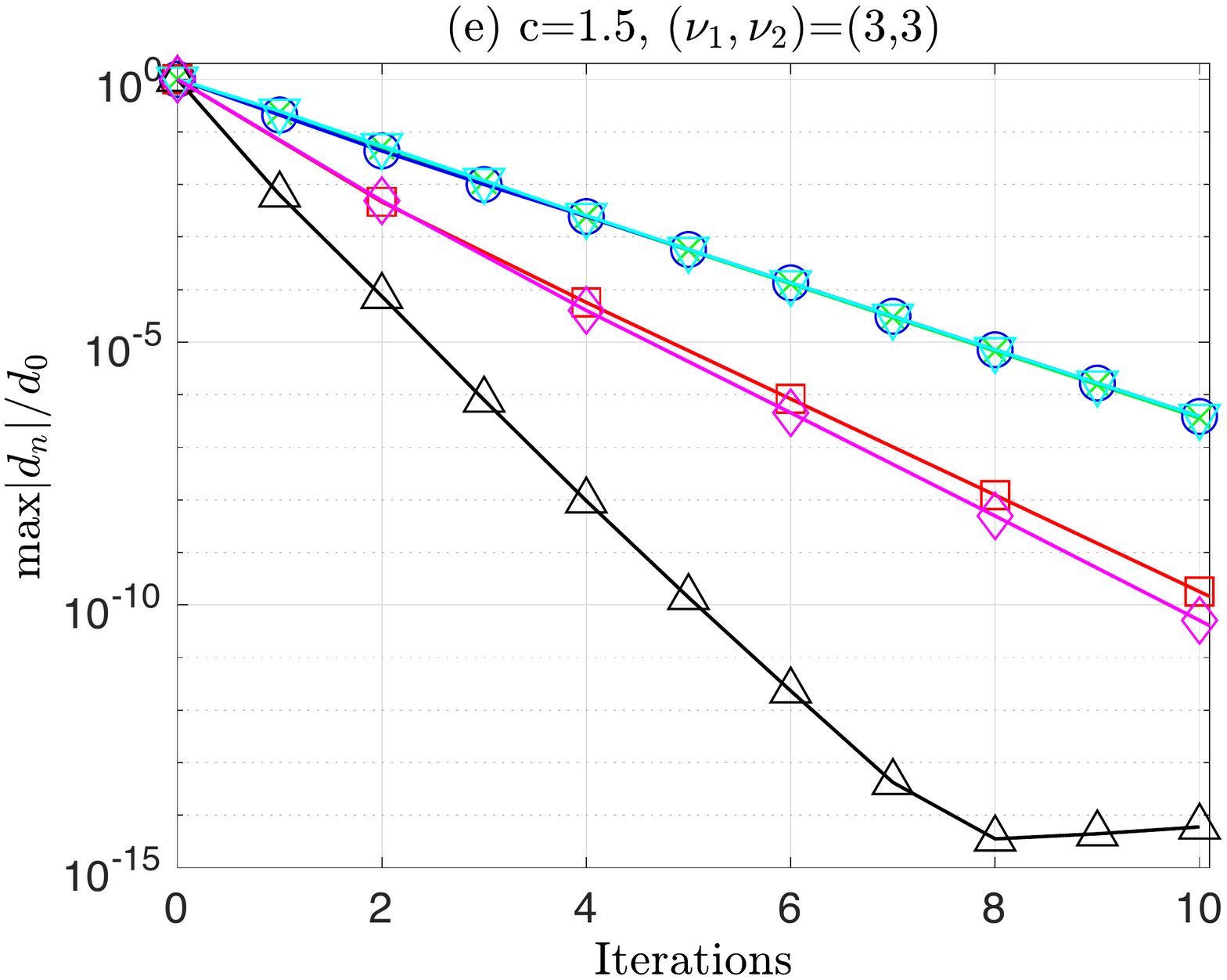}\\
\includegraphics[trim=3 0 40 0,clip,width=0.33\columnwidth]{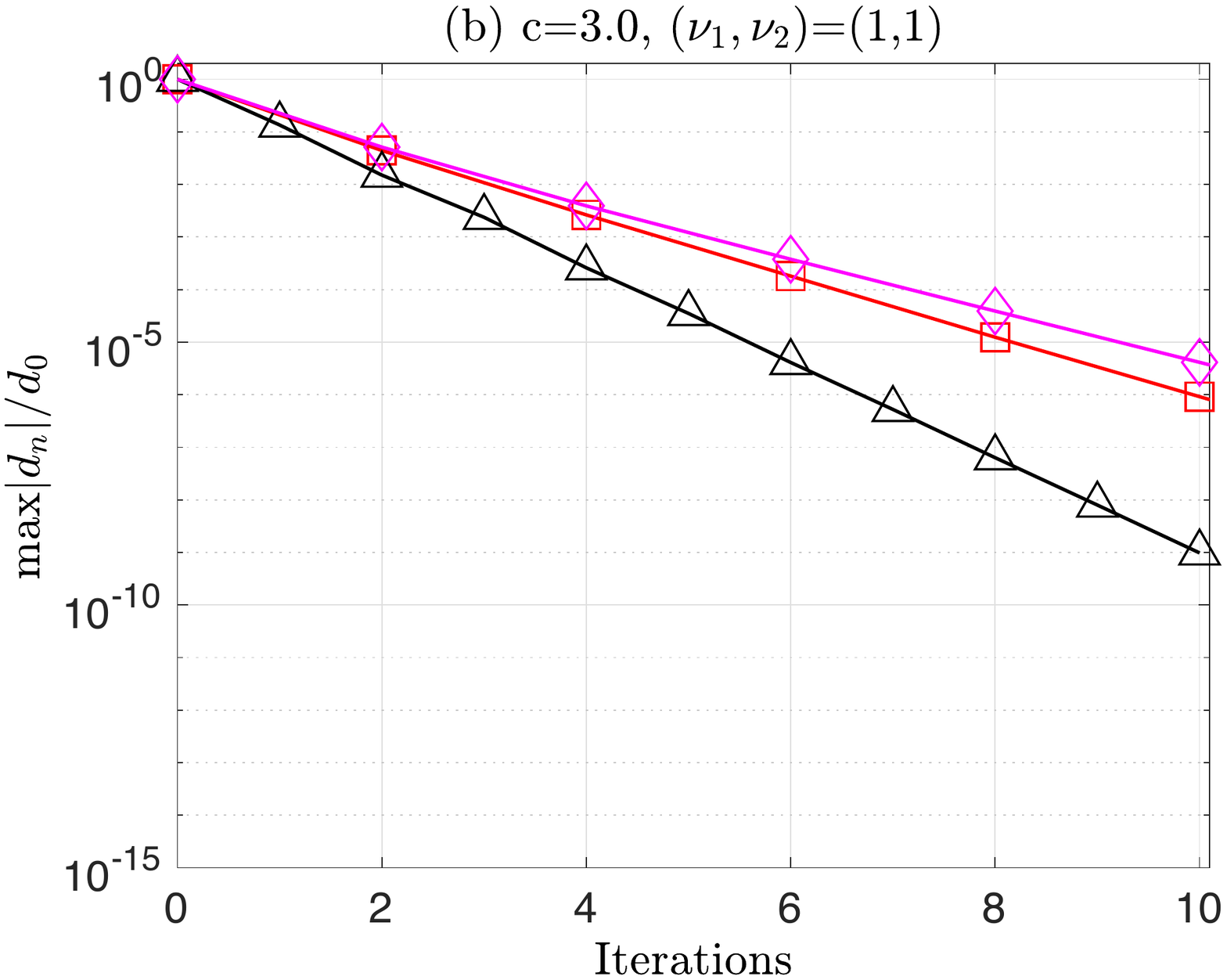}\includegraphics[trim=3 0 40 0,clip,width=0.33\columnwidth]{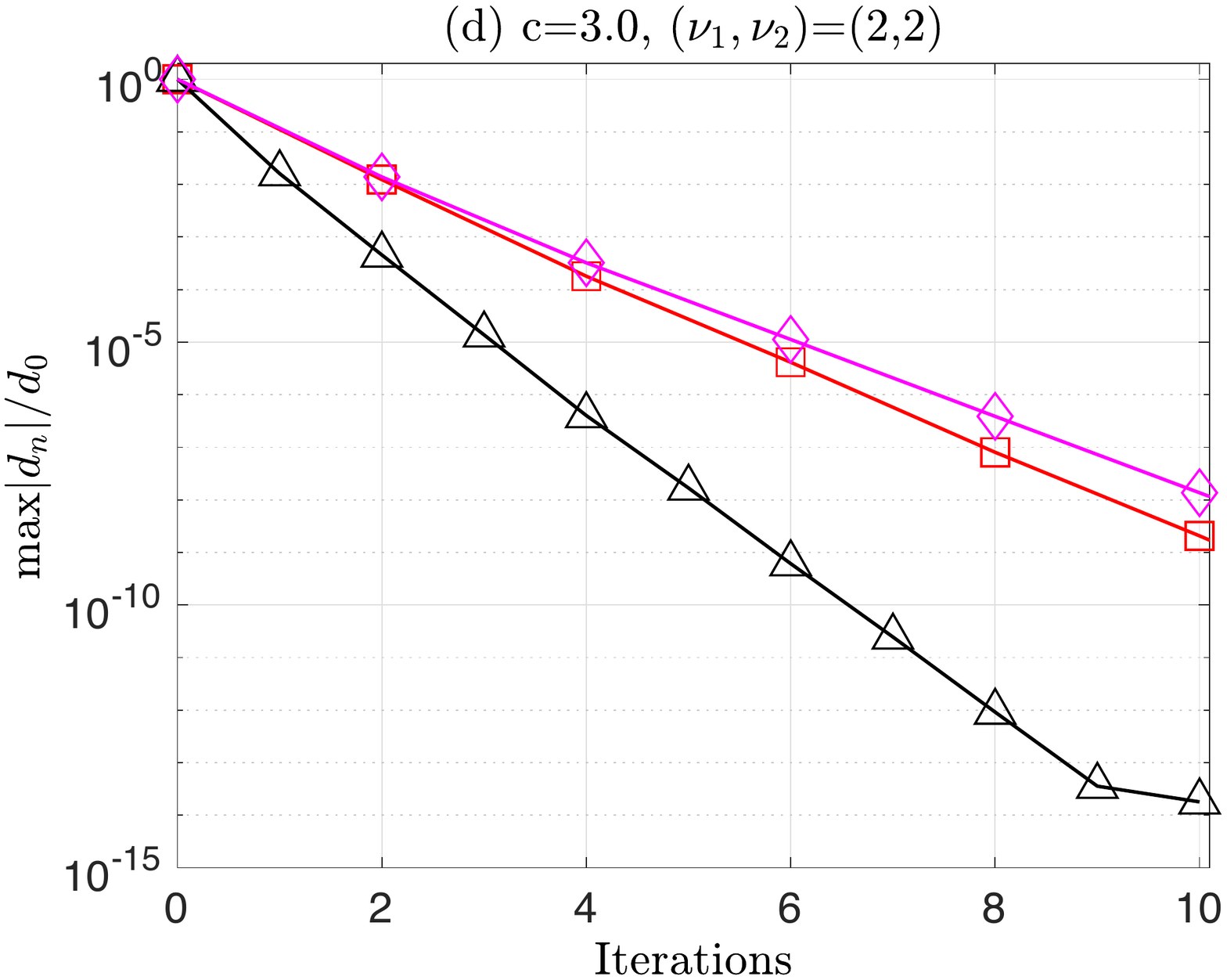}\includegraphics[trim=3 0 40 0,clip,width=0.33\columnwidth]{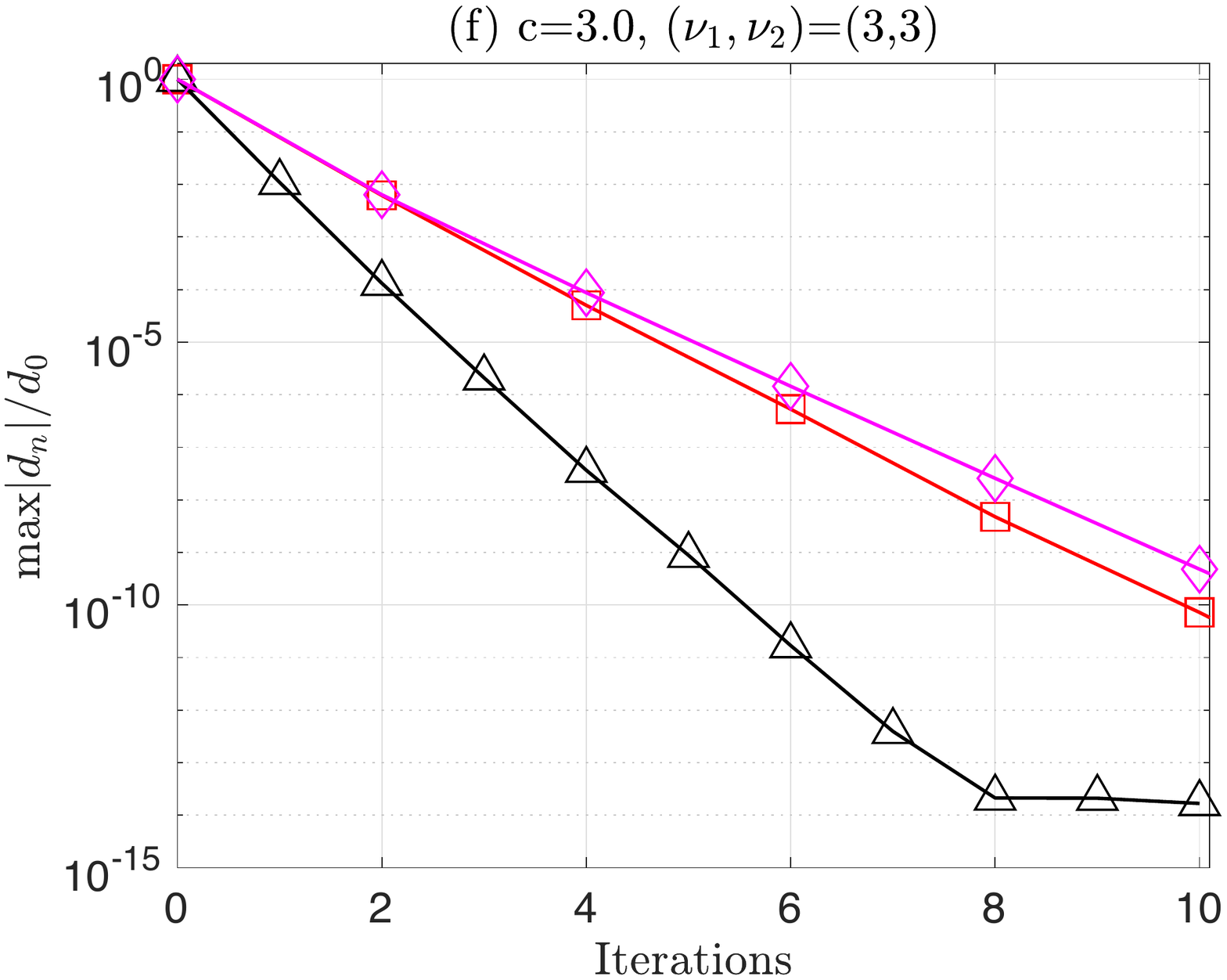}
\caption{Multigrid convergence on \eqref{prob:DP} over a $129 \times 129$ stretched grid with near-wall clustering and different smoothers:
(for legend, see Figure \ref{Fig8}).
(a,b,c) $c=1.5$; (d,e,f) $c=3.0$.
(a,d) $(\nu_1,\nu_2)=(1,1)$;
(b,e) $(\nu_1,\nu_2)=(2,2)$;
(c,f) $(\nu_1,\nu_2)=(3,3)$.\\
Note that checkerboard, one-direction zebra, and wireframe did not converge for the cases with $c=3.0$.  Tweed relaxation outperformed all other relaxation schemes tested.}
\label{Fig9}
\end{figure}

Tests on the 2D Poisson equation computed on a stretched grid with near-wall grid clustering (Figure \ref{Fig9}) demonstrate conclusively that,
per relaxation applied (which costs $\sim 12\,n^2$ flops in each case),
multigrid leveraging tweed relaxation is about twice as fast as
multigrid leveraging alternating-direction zebra relaxation, which is the next-best choice.  This result is achieved because, using the tweed motif, relaxation lines are locally orthogonal to the local directions of densest grid clustering {\it everywhere} in the computational domain.
Further, the amount of stretching applied appears to have a rather minor effect on the convergence rate obtained.

\begin{figure}  
\centering
\includegraphics[trim=3 0 40 0,clip,width=0.33\columnwidth]{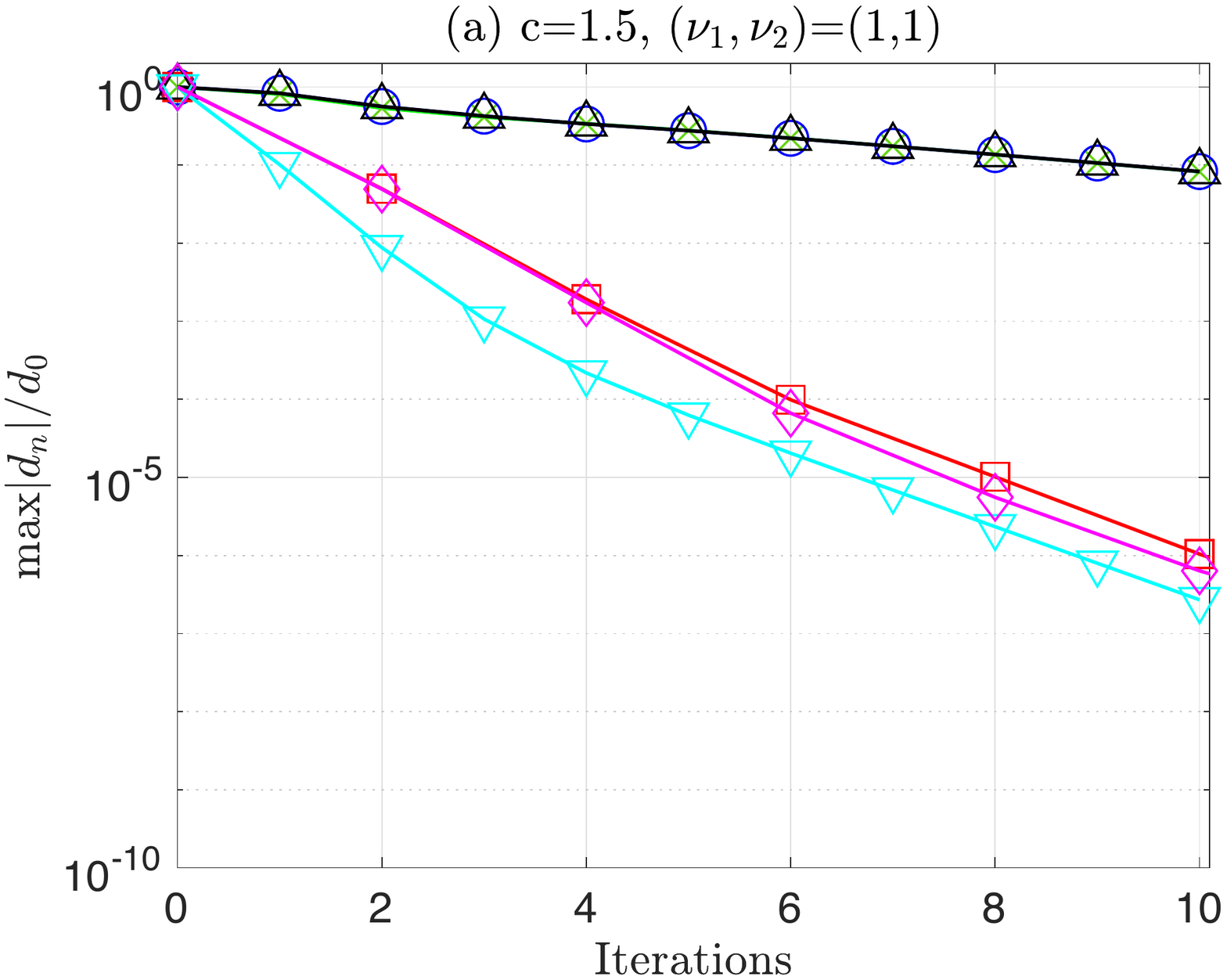}\includegraphics[trim=3 0 40 0,clip,width=0.33\columnwidth]{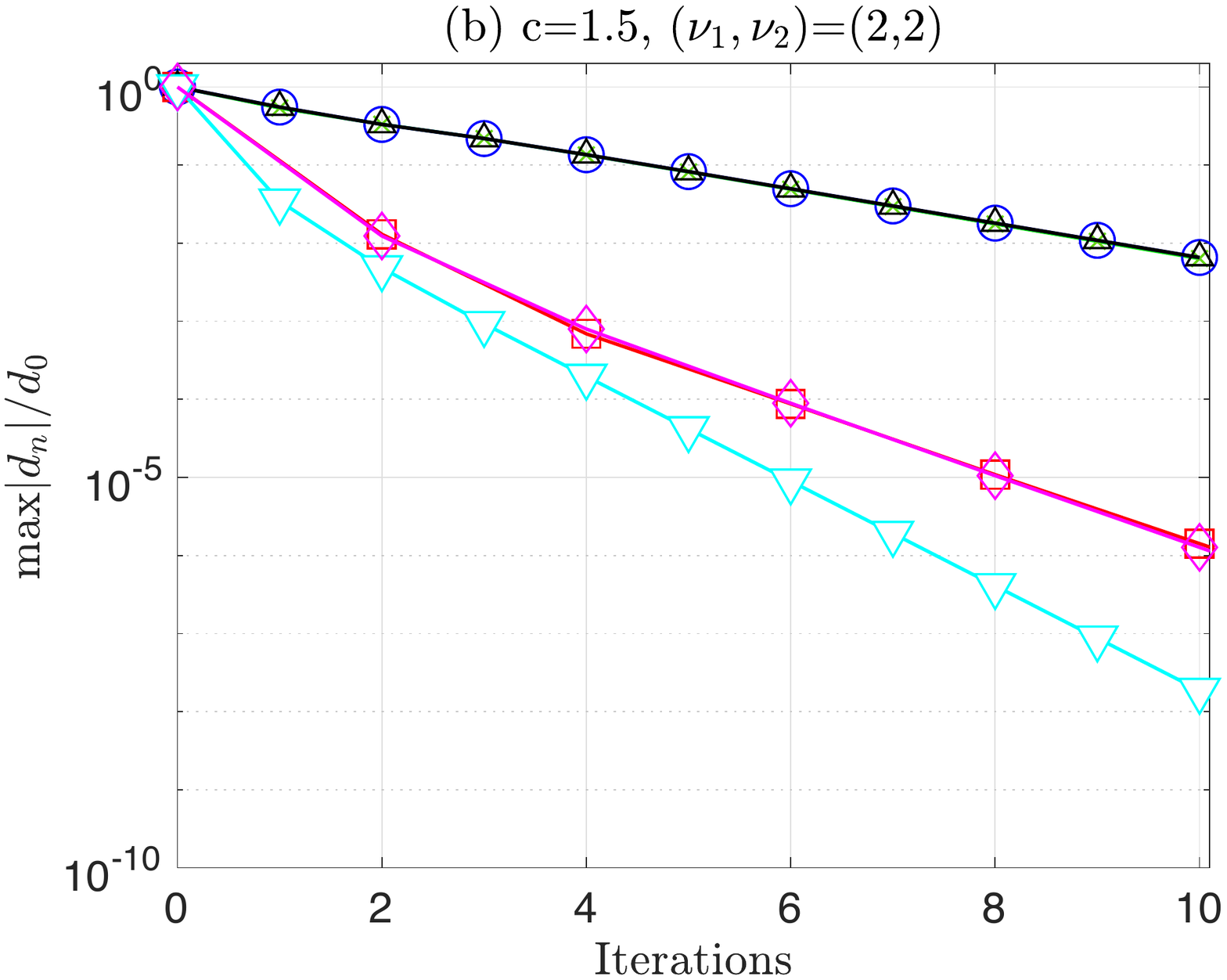}\includegraphics[trim=3 0 40 0,clip,width=0.33\columnwidth]{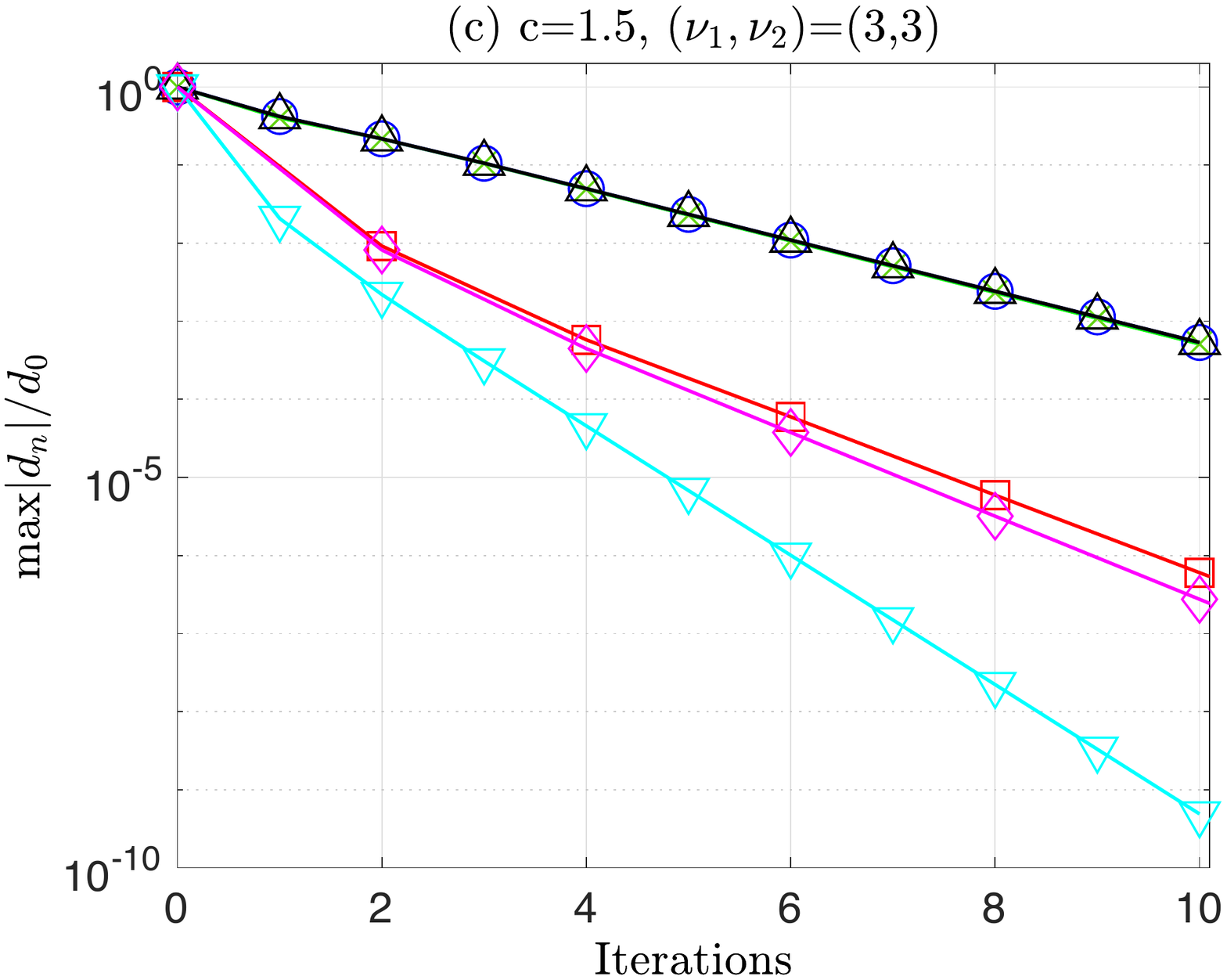}
\caption{Multigrid convergence on \eqref{prob:DP} over a $129 \times 129$ stretched grid with near-centre clustering and different smoothers
(for legend, see Figure \ref{Fig8}).
(a) $(\nu_1,\nu_2)=(1,1)$;
(b) $(\nu_1,\nu_2)=(2,2)$;
(c) $(\nu_1,\nu_2)=(3,3)$.  Wirefame relaxation slightly outperformed all other relaxation schemes tested.}
\label{Fig10}
\end{figure}

\begin{figure}  
	\centering
	\subfloat[\emph{Tweed}]{\psfrag{x}[t][b]{\small$\tilde{x}$}\psfrag{y}[b][c]{\small$\tilde{y}$}\includegraphics[width=0.45\columnwidth]{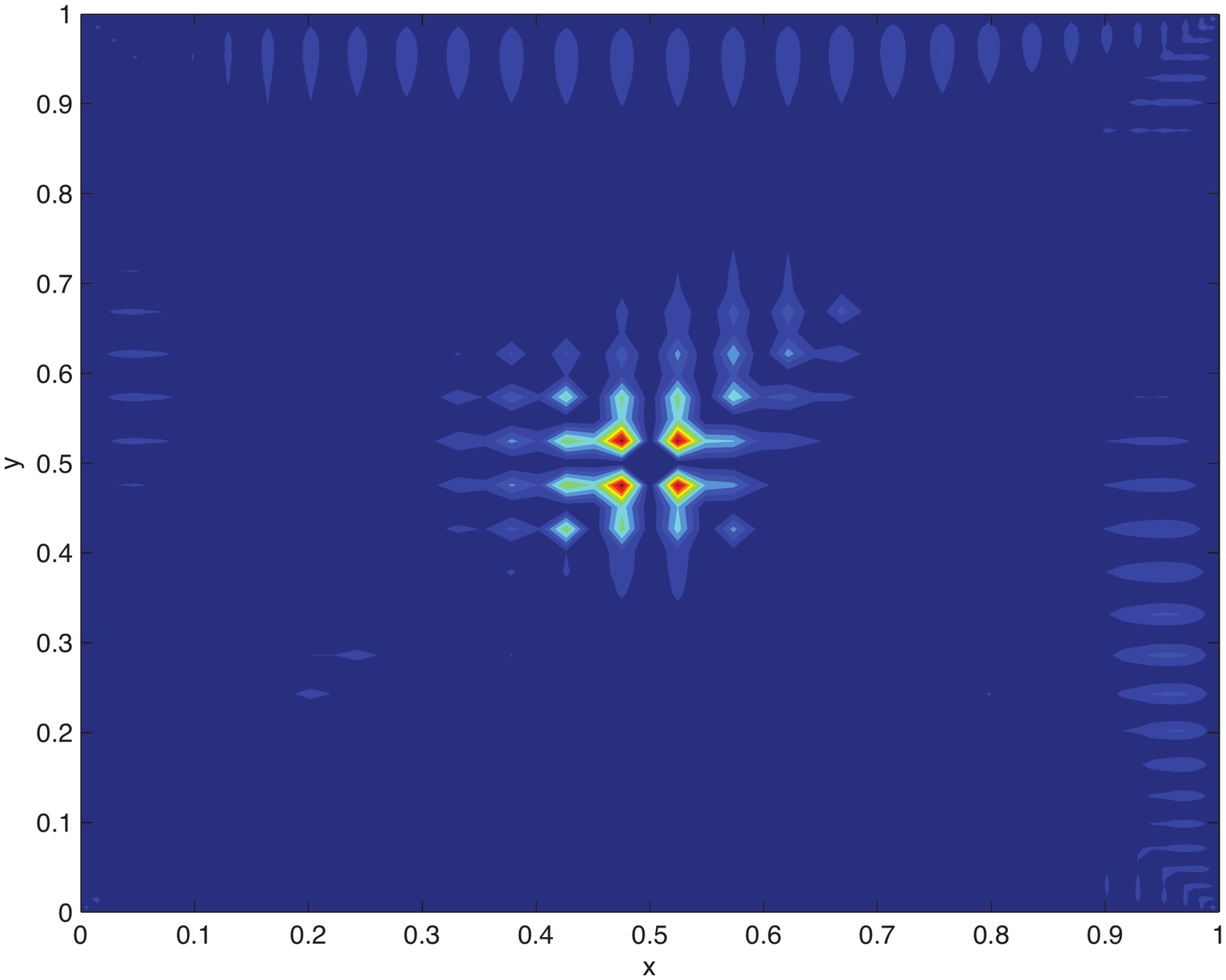}\label{fig:defect_tw}} \quad
	\subfloat[\emph{Alternating-direction zebra}]{\psfrag{x}[t][b]{\small$\tilde{x}$}\psfrag{y}[b][c]{\small$\tilde{y}$}\includegraphics[width=0.45\columnwidth]{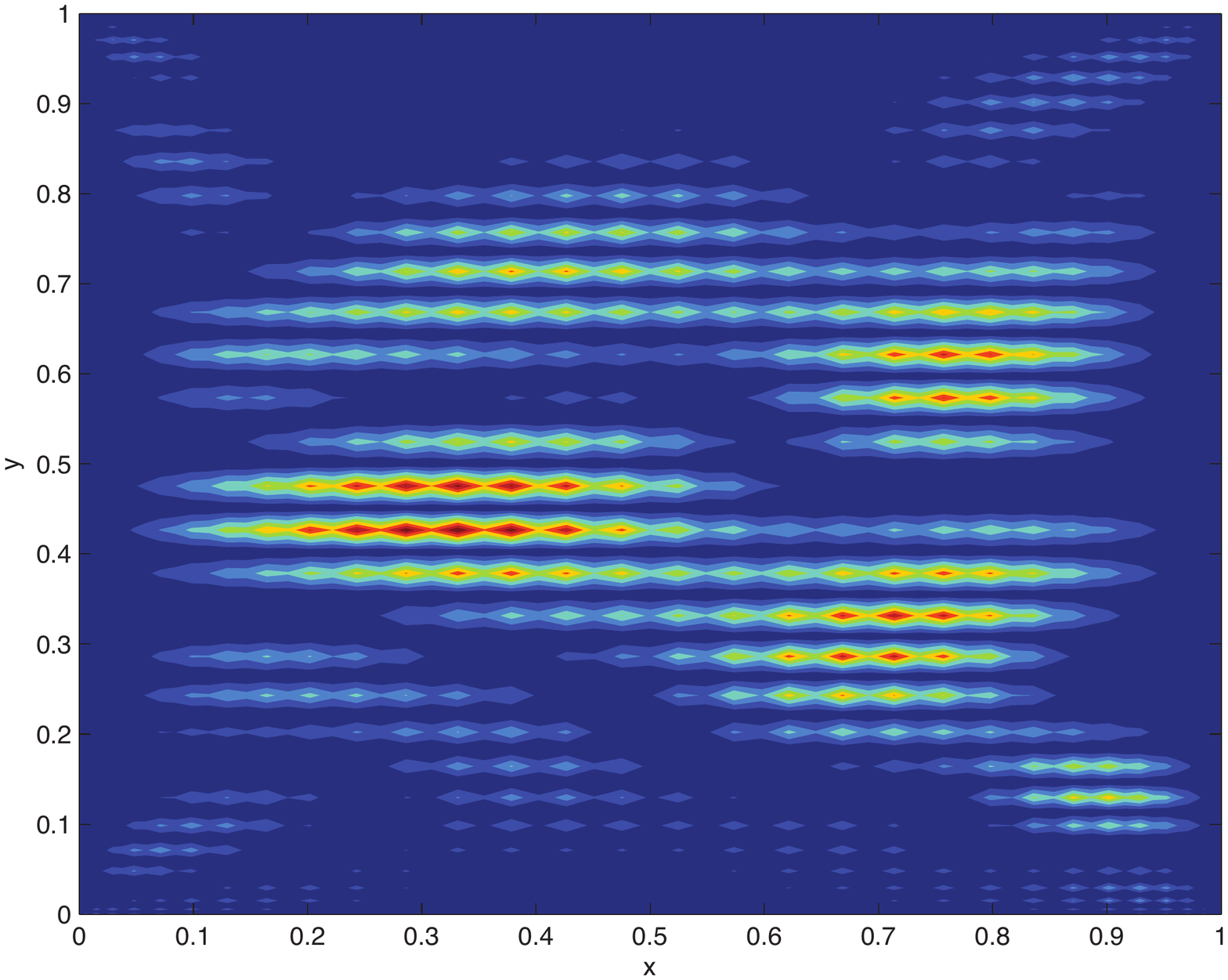}\label{fig:defect_zxy}}
	\caption{Defect after three cycles of multigrid over a $129 \times 129$ stretched grid with near-wall clustering [defined using~\eqref{eq:tanh}] with $(\nu_1,\,\nu_2) = (2,\,2)$ and different smoothing schemes.}
	\label{fig:defect}
\end{figure}

\begin{figure}  
	\centering
	\subfloat[\emph{Wireframe}]{\psfrag{x}[t][b]{\small$\tilde{x}$}\psfrag{y}[b][c]{\small$\tilde{y}$}\includegraphics[width=0.45\columnwidth]{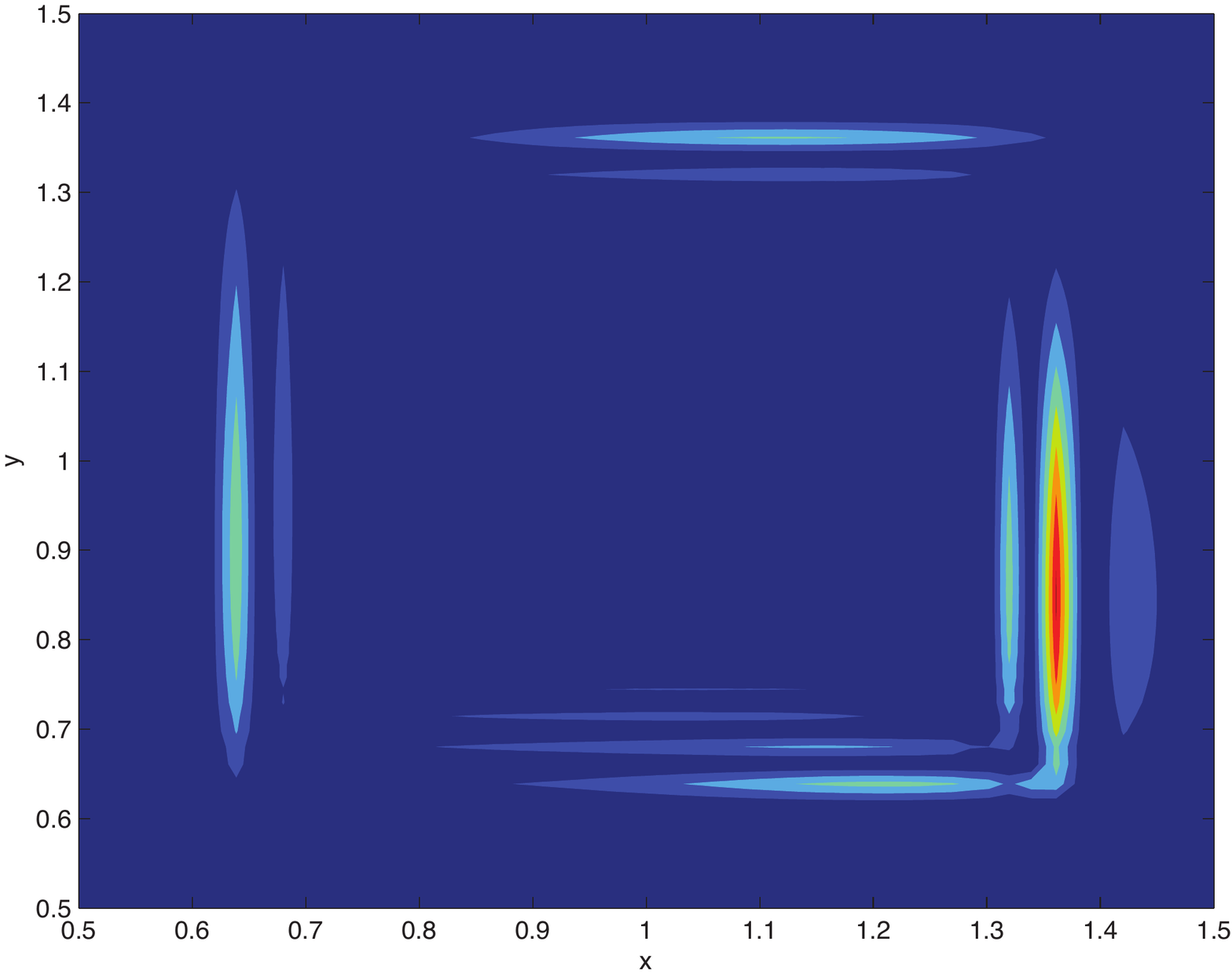}\label{fig:defectINV_box}} \quad
	\subfloat[\emph{Alternating-direction zebra}]{\psfrag{x}[t][b]{\small$\tilde{x}$}\psfrag{y}[b][c]{\small$\tilde{y}$}\includegraphics[width=0.45\columnwidth]{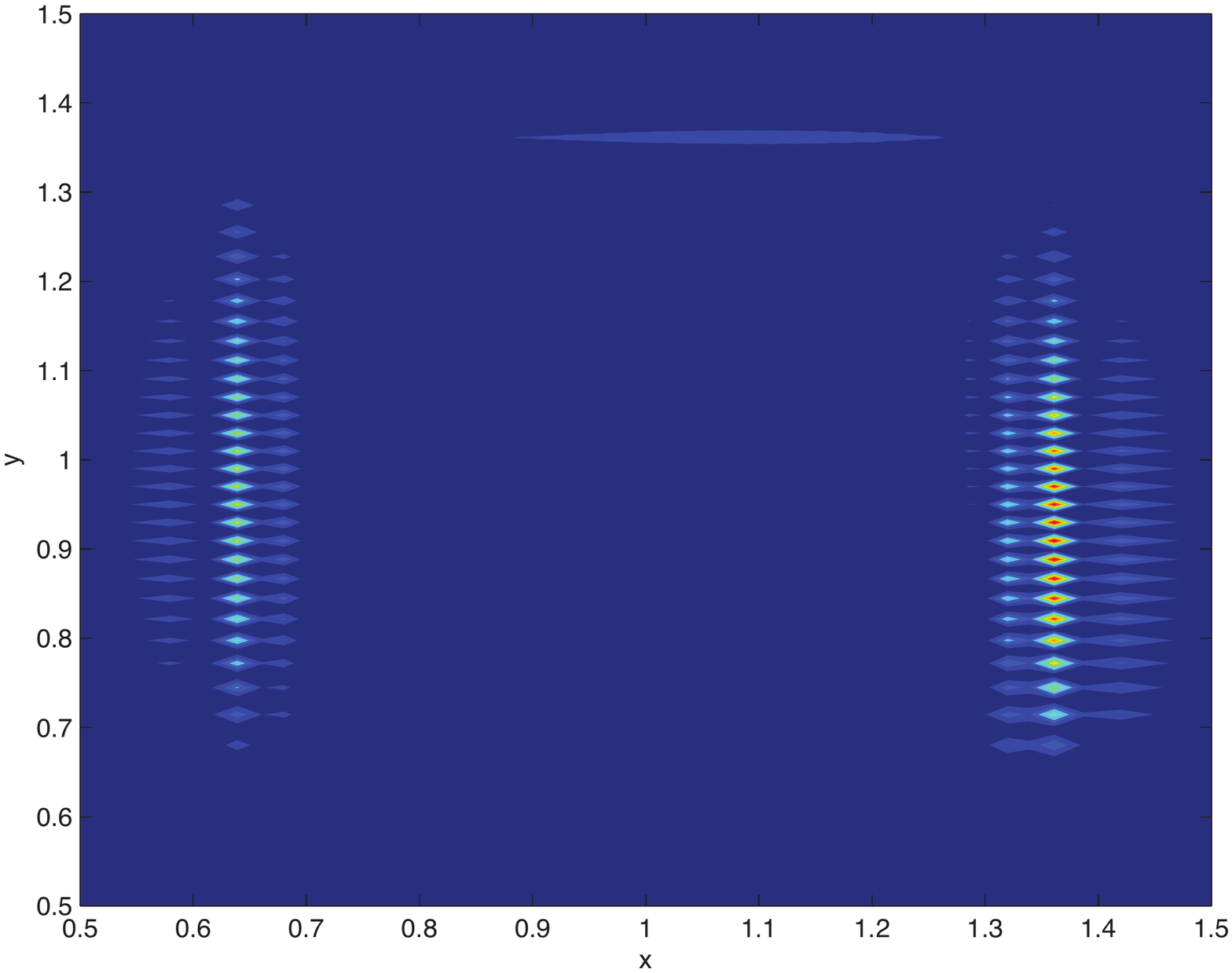}\label{fig:defectINV_zxy}}
	\caption{Defect after three cycles of multigrid over a $129 \times 129$ stretched grid with near-centre clustering [defined using~\eqref{eq:tanh_inv}] with $(\nu_1,\,\nu_2) = (2,\,2)$ and different smoothing schemes.}
	\label{fig:defectINV}
\end{figure}

\begin{figure}  
\centering
\includegraphics[width=0.5\columnwidth]{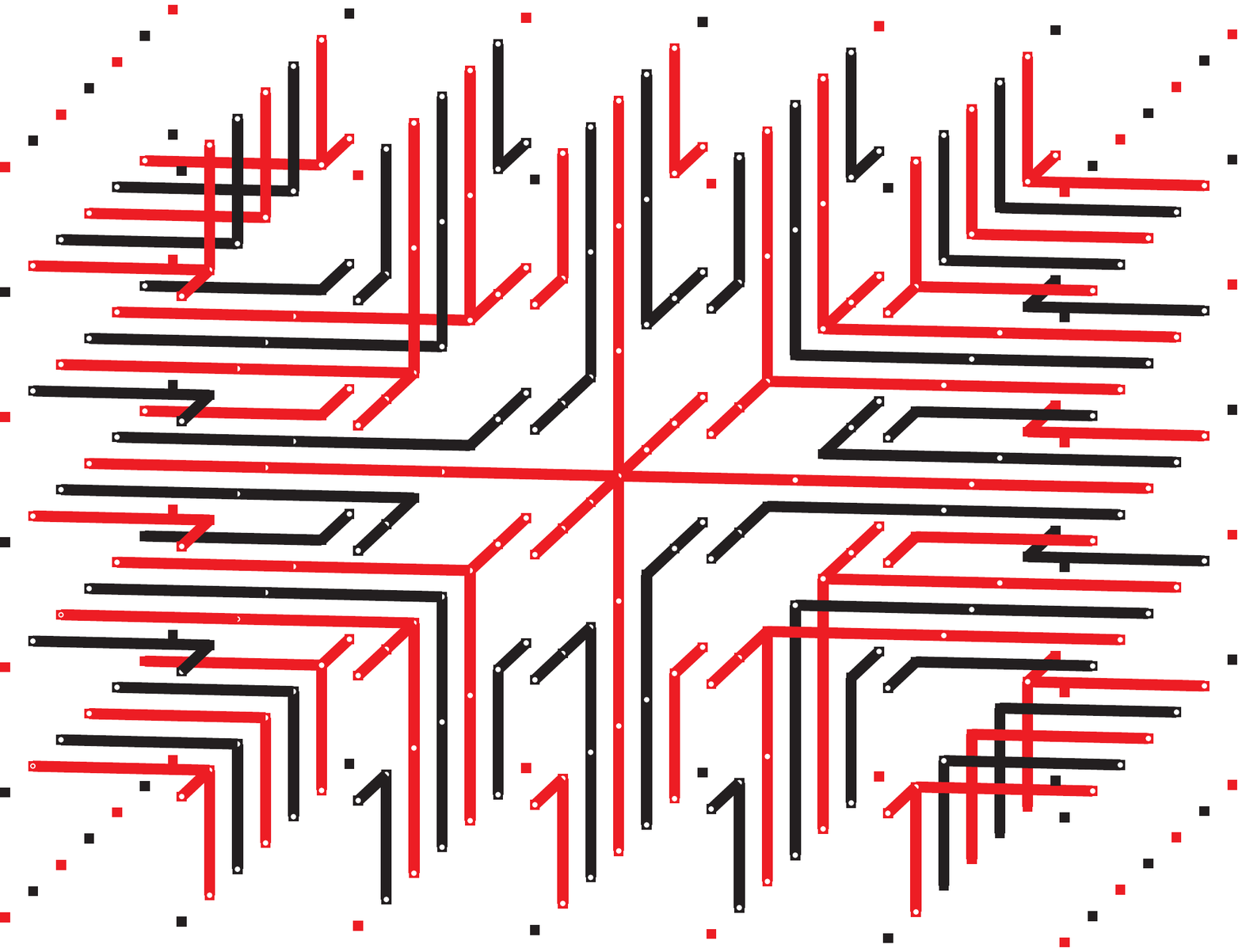}\\[0.01in]
\includegraphics[width=0.5\columnwidth]{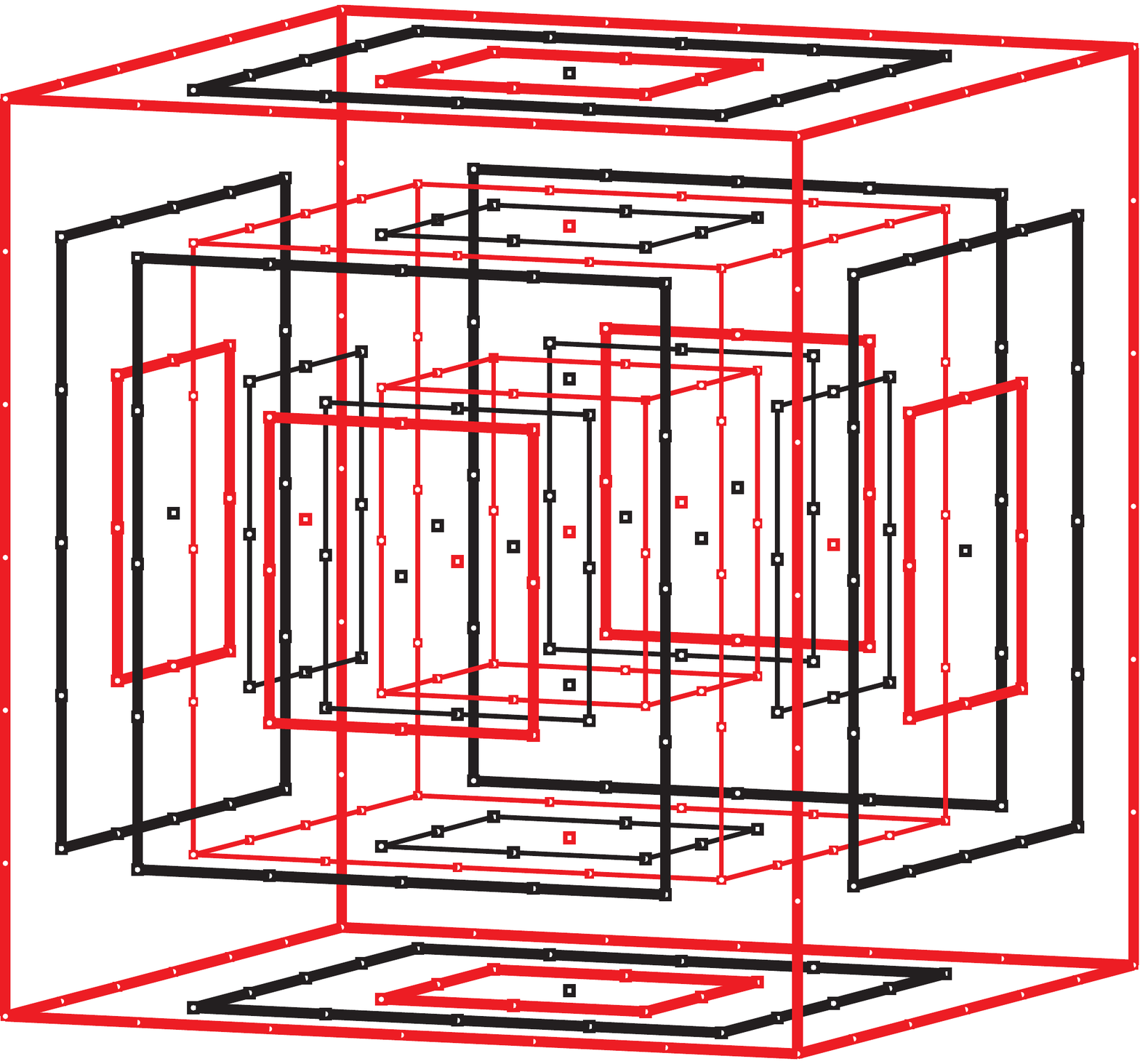}
\caption{Relaxation motifs for (top) tweed smoothing, and (bottom) wireframe smoothing, for 3D rectangular grids with $n_x = n_y = n_z$.  For the plot of the wireframe smoothing motif, to ease the visualization, the wireframes around the outermost cuboid are drawn with heavier lines than the wireframes around the cuboids on the interior.}
\label{Fig13}
\end{figure}

Tests on the 2D Poisson equation computed on a stretched grid with near-centre grid clustering (Figure \ref{Fig10}) demonstrate, more modestly, that
multigrid leveraging wireframe relaxation is at least competitive with
alternating-direction zebra relaxation.  Which is actually faster in practice will depend on machine-specific details, specifically the speed of floating-point computations versus the speed of memory references on the computer system used.

The extension of tweed and wireframe relaxation motifs to 3D is straightforward, as discussed in \S \ref{sec:3d}.  It has also been shown, in \S \ref{sec:Complementarity}, that the 2D tweed and wireframe relaxation motifs are in a sense complementary, and when superposed their blocks cover all of the connections in a 2D Cartesion grid.  This complementarity condition extends to 3D by defining a third relaxation motif, dubbed the {\it complementary wirefarme} motif, that agian covers the surface of concentric cuboids, with blocks that encircle the 8 corners of each cuboid rather than the 6 centers of each face.
This complementarity notion (in both 2D and 3D) adds a degree of geometric ``completeness'' to the entire tweed/wireframe framework, revealing that the motifs developed herein are perhaps, in a sense, a bit more fundamental than simply a convenient construction for a contrived numerical problem.

In addition to further quantifying the performance of the 3D extensions of the 2D tweed and 2D wireframe relaxation approaches (relative to 3D alternating-direction zebra, which is the next best choice), as described above, our team is currently implementing 2D and 3D tweed relaxation for the solution of the Poisson equation for the pressure update in the numerical solution of turbulent duct and cavity flows with grid clustering near the walls.  Results will be reported in future work.

%

\bibliographystyle{IEEEtran}
\bibliography{biblio}{}

\pagebreak
\newpage

\clearpage

\clearpage

\end{document}